
%
\def\PaperSize{letter}		
\input amstex.tex		
\documentstyle{amsppt}		
\def\TheMagstep{\magstep1}      
\def\papernumber{alg-geom/9412019}	

\def\TheABSTRACT{\ignorespaces
        We prove the results about mixed Buchsbaum--Rim multiplicities
announced in \cite{\KT, (9.10)(ii), \p.224}, including a general
mixed-multiplicity formula.  In addition, we identify these
multiplicities as the coefficients of the ``leading form'' of the
appropriate Buchsbaum-Rim polynomial in three variables, and we prove a
positivity theorem.  In fact, we define the multiplicities as the
degrees of certain zero-dimensional ``mixed twisted'' Segre classes,
and we develop an encompassing general theory of these new rational
equivalence classes in all dimensions.  In parallel, we develop a
theory of pure ``twisted'' Segre classes, and we recover the main
results in \cite{\KT} about the pure Buchsbaum--Rim multiplicities, the
polar multiplicities, and so forth.  Moreover, we identify the
additivity theorem \cite{\KT, (6.7b)(i), \p.205} as giving a sort of
residual-intersection formula, and we show its (somewhat unexpected)
connection to the mixed-multiplicity formula.  Also, we work in a more
general setup than in \cite{\KT}, and we develop a new approach, based
on the completed normal cone.
  }

\def\RefKeys{AM BR Fulton FL EGA KT NR Teissier Thorup } 
\def\TRUE{TRUE}
\ifx\DoublepageOutput\TRUE \def\TheMagstep{\magstep0} \fi
\mag=\TheMagstep
\TagsOnRight

\abovedisplayskip 6pt plus6pt minus0.25pt
\belowdisplayskip=\abovedisplayskip
\abovedisplayshortskip=0mm
\belowdisplayshortskip=2mm
\def\centertext
 {\hoffset=\pgwidth \advance\hoffset-\hsize
  \advance\hoffset-2truein \divide\hoffset by 2
  \voffset=\pgheight \advance\voffset-\vsize
  \advance\voffset-2truein \divide\voffset by 2
 }
\newdimen\pgwidth	\newdimen\pgheight
\def\letter{letter}	\def\AFour{AFour}
\ifx\PaperSize\letter
	\pgwidth=8.5truein \pgheight=11truein
 \message{- Got a papersize of letter.  }\centertext \fi
\ifx\PaperSize\AFour
 \pgwidth=210truemm \pgheight=297truemm
 \message{- Got a papersize of AFour.  }\centertext \fi

 \newdimen\fullhsize 	\newbox\leftcolumn
 \def\fulline{\hbox to \fullhsize}
\def\doublepageoutput
{\let\lr=L
 \output={\if L\lr
          \global\setbox\leftcolumn=\columnbox \global\let\lr=R%
        \else \doubleformat \global\let\lr=L\fi
        \ifnum\outputpenalty>-20000 \else\dosupereject\fi}%
 \def\doubleformat{\shipout\vbox{%
        \fulline{\hfil\box\leftcolumn\hfil\columnbox\hfil}%
	}%
 }%
 \def\columnbox{\vbox
   {\makeheadline\pagebody\makefootline\advancepageno}%
   }
 \fullhsize=\pgheight \hoffset=-1truein
 \voffset=\pgwidth \advance\voffset-\vsize
  \advance\voffset-2truein \divide\voffset by 2
\headline={%
 \eightpoint
  \ifnum\pageno=1\firstheadline
  \else
    \ifodd\pageno\rightheadline
    \else\leftheadline\fi
  \fi
}
\let\firstheadline=\hfil

\null\vfill\nopagenumbers\eject\pageno=1\relax 
}
\ifx\DoublepageOutput\TRUE \doublepageoutput \fi

\def\today{\number\day \space\ifcase\month\or
 January\or February\or March\or April\or May\or June\or
 July\or August\or September\or October\or November\or December\fi
\space \number\year}

 \leftheadtext{KLEIMAN AND THORUP}
 \rightheadtext{MIXED BUCHSBAUM--RIM MULTIPLICITIES}
 \def\rightheadline{\rlap{\papernumber}\hfil
	\the\rightheadtoks\hfil\llap{\folio}}
 \def\leftheadline{\rlap{\folio}\hfil
 	\the\leftheadtoks\hfil\llap{\papernumber}}

 \newskip\sectskipamount \sectskipamount=0pt plus30pt
 \def\sectionhead#1 #2\par{\def\sectno{#1}\def\sectname{#2}%
   \vskip\sectskipamount\penalty-250\vskip-\sectskipamount
   \bigskip
   \centerline{\smc \number\sectno.\enspace\sectname}\nobreak
   \medskip
   \message{\number\sectno. \sectname }%
}

 \def\proclaimm#1#2 {\medbreak\noindent {\bf(\sectno.#2) #1.}\quad
	\begingroup\it}
 \def\subh#1#2{\medbreak\noindent {\bf(\sectno.#2) #1.}\quad}
\def\art#1 #2\par{\subh{\rm({\it#2\unskip\/})}{#1}}

\def\stp{\subh{Setup}}
\def\expl{\subh{Example}}
\def\prop{\proclaimm{Proposition}}
	\def\propx #1{\proclaimm{Proposition {\rm (#1)}}}
\def\thm{\proclaimm{Theorem}}
	\def\thmx #1{\proclaimm{Theorem {\rm (#1)}}}
\def\cor{\proclaimm{Corollary}}
\def\lem{\proclaimm{Lemma}}
	\def\lemx #1{\proclaimm{Lemma {\rm (#1)}}}
\def\pf{\endgroup\medbreak\noindent {\bf Proof.}\quad}
\def\demobox{\vbox{\hrule\hbox{\vrule\kern.5ex
	\vbox{\kern1.2ex}\vrule}\hrule}}
\def\enddemo{{\unskip\nobreak\hfil\penalty50
  \hskip1em\hbox{}\nobreak\hfil\demobox
  \parfillskip=0pt \finalhyphendemerits=0 \par}}

\def\Cs#1){{\rm(\sectno.#1)}} \def\dfn#1{{\it #1}}
\def\and{\hbox{ and }}	
\def\tgs{\tag\sectno.}
\def\part #1 {\par{\rm(#1)~~}}

\let\wh=\widehat	\let\wt=\widetilde
\let\onto=\twoheadrightarrow
\let\into=\hookrightarrow
\let\scq=\succcurlyeq
\def\risom{\,\overset{\sim}\to{\smash\longrightarrow\vrule
    height0.5ex width0pt depth0pt}\,}
\let\ox=\otimes	\let\?=\overline
\def\uprab#1{^{\langle#1\rangle}}

\def\opdef#1 {\expandafter\def\csname#1\endcsname{\operatorname{#1}}}
\opdef Grass	\opdef rk	\opdef td	\opdef Cyc
\opdef Spec	\opdef Supp	\opdef Proj	\opdef length
\opdef ord	\opdef Ker	\opdef Cok	\opdef cod
\opdef Im	\opdef A

\def\cG{{\Cal G}}	\def\K{{\Cal K}}	\def\O{{\Cal O}}
\def\L{{\Cal L}}	\def\M{{\Cal M}}	\def\N{{\Cal N}}
\def\I{{\Cal I}}	\def\J{{\Cal J}}

		\def\S{\bold S}		\def\T{\bold T}
\def\IP{\bold P}	

\def\atref#1{\ref\key\csname#1\endcsname}
 \newcount\refno \refno=0
 \def\MakeKey{\advance\refno by 1 \expandafter\xdef
 	\csname\TheKey\endcsname{{%
		\number\refno}}\NextKey}
 \def\NextKey#1 {\def\TheKey{#1}\ifx\TheKey\NoKey\let\next\relax
  \else\let\next\MakeKey \fi \next}
 \def\NoKey{*!*}
 \expandafter\NextKey \RefKeys *!*
 \ifx\UseNumericalRefKeys\TRUE \widestnumber\no{\number\refno}\fi
\def\p.{\unskip\space p.\penalty10000 \thinspace}
\def\pp.{\unskip\space pp.\penalty10000 \thinspace}


\topmatter
 \title
	Mixed Buchsbaum--Rim Multiplicities
 \endtitle
 \author
	Steven Kleiman$^*$ and Anders Thorup$^\dagger$
 \endauthor
 \address
   Department of Mathematics, Massachusetts Institute of Technology,
   Cambridge, MA 02139, USA
 \endaddress
 \email \tt kleiman\@math.mit.edu \endemail
 \address
   Matematisk Institut, Universitetsparken 5,
   DK-2100 Copenhagen, Denmark
 \endaddress
 \email \tt thorup\@math.ku.dk \endemail
 \date\papernumber\enddate
 \thanks
	$^*$Supported in part by NSF Grants 9106444-DMS and
9400918-DMS \endgraf
	$^\dagger$Supported in part by Danish Natural Science Research
	Council Grant 11-7428 \endthanks
 \abstract\TheABSTRACT\endabstract
\endtopmatter

\document

\sectionhead 1 Introduction

The theory of mixed multiplicities of primary ideals was introduced by
Teissier in his study of complex hypersurface germs with isolated
singularities.  A decade later, Gaffney began extending Teissier's work
to complete intersections, and was led to conjecture a theory of
generalized multiplicities of submodules of finite colength in a free
module, including an important mixed-multiplicity formula for the
product of an ideal and a submodule.  It turned out that these
generalized multiplicities are nothing but the multiplicities introduced
a decade before Teissier's work by Buchsbaum and Rim, who established
many of their fundamental properties, but no mixed-multiplicity
formula.  Recently, the authors gave a general treatment of the
Buchsbaum--Rim multiplicity, based on blowups and intersection numbers,
in \cite{\KT} (that paper also contains a more extensive history of the
subject).  On \p.225, the authors announced a mixed-multiplicity
formula for an arbitrary pair of submodules.  Here we prove an even
more general formula, and show that it's closely related, surprisingly,
to another fundamental formula, the additivity formula.  We also
simplify, generalize, and advance the previous treatment via a
new approach.

Sections~2 and 3 study two preliminary notions, module transforms and
distinguished subsets.  Section~4 studies a ``twisted'' version of the
usual Segre classes of a subscheme.  The degrees of these classes yield
the Buchsbaum--Rim multiplicities.  Moreover, the usual theory of Segre
classes leads to simpler and more natural proofs of their properties.
One case in point is the additivity theorem (4.6), one of our main
results.  It yields the old additivity formula for Buchsbaum--Rim
multiplicities \cite{\KT, (6.7b)(i), \p.205}, which had a somewhat
mysterious intersection-theoretic proof.  By contrast, the
intersection-theoretic proof of (4.6) is inspired by the proof of the
residual intersection theorem \cite{\FL, Prop.~2, \p.174}.  Section~5
studies formally similar ``mixed twisted'' Segre classes, which
generalize the mixed Buchsbaum--Rim multiplicities.  Theorem~(5.5),
another main result, gives the generalization of the announced
mixed-multiplicity formula. The expansion formula (5.6.1) expresses the
mixed classes in terms of the pure classes; given it, we show in (5.7)
that the mixed-class formula and the additivity formula are, somewhat
surprisingly, essentially equivalent.  Section~6 determines the result
of intersecting with a suitably general pseudo-divisor in preparation
for Section~7.  Section~7 contains the deepest results, which concern
positivity.  The results in Sections~6 and 7 are modeled on results in
\cite{\KT}, which they generalize and improve.  Finally, Section~8
shows that the Buchsbaum--Rim multiplicities, both pure and mixed, are
the leading coefficients of generalized Hilbert polynomials, which we
call {\it Buchsbaum--Rim polynomials.}

We now discuss the contents in more detail.  We wish, first of all, to
survey the paper, and describe the main constructions, results, and
proofs.  We wish, at the same time, to compare the new treatment with
the old, and highlight the advantages of the new approach.  For the
most part, we save until (7.10) and (8.5) the matter of how to formally
recover the old results from the new ones.  The two treatments are
logically independent, although we do occasionally refer to the old one
for some preliminary general lemmas of dimension theory and the theory
of Hilbert polynomials.  In both treatments, modern intersection theory
is an essential tool.  Since we work on arbitrary Noetherian schemes,
the development in Fulton's book \cite{\Fulton} is not sufficiently
general, but a suitable extension is found in \cite{\Thorup}.

To introduce the notion of Buchsbaum-Rim multiplicities and our new
approach to the theory, it is best to begin the discussion with the
last section, Section~8, and to dwell on it longer than we would
otherwise.  Let then $A$ be an arbitrary Noetherian local ring, $G$ a
graded $A$-algebra that is generated by finitely many elements of
degree one, $H$ an $A$-submodule of $G_1$, and $M$ a finitely generated
graded $G$-module.   As a function of $p$ and $n$, the length,
	$$\lambda(p,n):=\length(M_{p+n}/H^pM_n),$$
 is eventually a polynomial, and its total degree is at most the
dimension $r$ of the sheaf $\M:=\wt M$ on $X:=\Proj(G)$, provided the
closed fiber of $X$ contains the intersection $Z\cap\Supp\M$, where $Z$
is the subscheme of $X$ defined by $H$.  This theorem was proved in
\cite{\KT, (5.10), \p.199}, and it is recovered in (8.2).  (Beware,
however, the notation has changed: in \cite{\KT}, we denoted $\Proj(G)$
by $P$ and $\Spec(A)$ by $X$.)  In fact, (8.2) proves a more general
version, in which $H\subseteq G_d$ for an arbitrary $d\ge0$.
Curiously, the case $d=0$ holds special interest, because it recovers a
result about the polar multiplicities, discussed in the middle of
\cite{\KT, (8.5), \p.214}.

The proof of (8.2) is significantly shorter than that of \cite{\KT,
(5.10)}.  Indeed, $M_{p+n}/H^{p+1}M_{n-1}$ has a natural filtration,
whose factors are the $A$-modules,
      $$N_{\nu,\mu}:=H^\nu M_\mu/H^{\nu+1}M_{\mu-1},$$
 for $\nu=0,\dots,p$ and $\mu=n+p-\nu$.  In \cite{\KT}, we worked with
the bigraded module whose $(p,q)$-th piece is $N_{p,q}$.  In this
paper, we work with the bigraded module $P_H(M)$ whose $(p,q)$-th piece
$P_H(M)_{p,q}$ is the direct sum $\bigoplus_{\nu=0}^p N_{\nu,p-\nu+q}$.
Here $P_H(M)$ is a module over $P_H(G)$, a natural bigraded
$A$-algebra.  Hence $P_H(M)$ gives rise to a sheaf on the scheme
arising from $P_H(G)$; in fact, this scheme is the completed normal
cone $P_Z$ of $Z$ in $X$, and so we denote the sheaf by $P_Z(\M)$.
Hence, by the theory of Hilbert polynomials, the length of
$P_H(M)_{p,q}$ is eventually a polynomial of total degree at most $r$;
whence, so is $\lambda(p,n)$.  In \cite{\KT}, instead of $P_Z$, we used
the projectivized normal cone, viewed as the exceptional divisor $D$ of
the blowup of $X$ along $Z$, and the proof was more than twice as long!

The proof of (8.2) yields more: the function $\lambda(p,n)$ eventually
begins with the binary form,
       $$\Lambda(p,n):=\sum_{i+k=r}e^{i,k}(\M) p^i n^k/i!\,k!,$$
 whose coefficients $e^{i,k}(\M)$ are equal to certain intersection
numbers; namely,
       $$e^{i,k}(\M)=\int c_1\O_X(1)^k c_1\O_{P_Z}(1)^i[P_Z(\M)]_r.$$
 On the other hand, these coefficients are, by \cite{\KT, (7.5),
\p.210}, equal to the associated multiplicities of $\M$; in
particular, $e^{r,0}(\M)$ is equal to the Buchsbaum--Rim multiplicity.
In \cite{\KT, (7.1), \p.207}, these multiplicities were defined as
sums of certain Segre numbers, which are intersection numbers on $X$
and $D$.  Hence these sums of intersection numbers are equal to single
intersection numbers, namely, those displayed.  These equalities also
have a purely intersection-theoretic proof; see (7.10).

There is a ``mixed'' version of the above theory.  Consider two closed
subschemes $Z_1$ and $Z_2$ of $X$ defined by $A$-submodules $H_1$ of
$G_{d_1}$ and $H_2$ of $G_{d_2}$ with $d_1,d_2\ge0$.  Then, as a
function of $p$, $q$ and $n$, the length,
	$$\lambda(p,q,n):=\length(M_{d_1p+d_2q+n}/H^{p}_1H^{q}_2M_n),$$
 is eventually a polynomial of total degree at most $r$, provided the
closed fiber of $X$ contains the intersection $Z\cap\Supp\M$.  The term
of total degree $r$ is the form,
	$$\Lambda(p,q,n):=\sum_{i+j+k=r}
		e^{i,j,k}(\M) p^i q^j n^k/i!\,j!\,k!,$$
 whose coefficients are equal to the mixed Buchsbaum--Rim
multiplicities of $\M$; these multiplicities are the intersection
numbers displayed in the next paragraph.  This theorem is (8.4).
Again, the case $d_1=0$ holds special interest because it recovers the
main conclusions of \cite{\KT, (9.10)(i), \pp.223--4}, which dealt
with an ideal and a submodule.  For arbitrary $d_1$, the numbers
$e^{i,j,k}(\M)$ are mixed versions of the associated multiplicities.
In particular, $e^{i,j,0}(\M)$ is equal to the mixed multiplicity
$e^{i,j}(\S)$ of \cite{\KT, (9.10)(ii), \p.224} with $\S:=[\M]_r$; see
(7.10).

The proof of (8.4) is analogous to that of (8.2).  Namely,
	$$M_{d_1p+d_2q+n}/H^{p+1}_1H^{q+1}_2M_{n-d_1-d_2}$$
 has a certain filtration, and the direct sum of its factors is the
$(p,q,n)$-th graded piece of a suitable trigraded module.  However, the
filtration is more involved this time: it is defined by taking sums of
products of the ideals of $Z_1$ and $Z_2$; see (2.7).  The trigraded
module gives rise to a sheaf $P^{1,2}_Z(\M)$ on a scheme $P^{1,2}_Z$.
This scheme is a mixed version of $P_Z$ where $Z:=Z_1+Z_2$.  Finally,
the theory of Hilbert polynomials yields the result with
	$$e^{i,j,k}(\M):=\int c_1\O_X(1)^k
   c_1\O_{P^{1,2}_Z}(1,0)^ic_1\O_{P^{1,2}_Z}(0,1)^j[P^{1,2}_Z(\M)]_r$$
 where $\O_{P^{1,2}_Z}(1,0)$ and $\O_{P^{1,2}_Z}(0,1)$ are the two
tautological invertible sheaves.

Thus there are two good reasons for using the completed cone $P_Z$ and
its mixed form $P^{1,2}_Z$ instead of the projectivized cone $D$ and
its mixed form $D_{1,2}$.  First, the treatment of the Buchsbaum--Rim
polynomials is shorter, simpler and more natural.  Second, the
Buchsbaum--Rim multiplicities are given by single intersection numbers,
rather than by unwieldy sums.  What is more, the standard theory of
Segre classes \cite{\Fulton, 4.2, pp. 73--79} now uses $P_Z$ instead of
$D$.

In Sections 2 to 7, the setup is more general: $X$ is an arbitrary
Noetherian scheme of finite dimension, $Z$ is a closed subscheme of
$X$, and $\M$ is a coherent $\O_X$-module.  In this generality, the
transforms $P_Z(\M)$ and $P^{1,2}_Z(\M)$ are constructed and studied in
Section~2 using a third and more fundamental transform $B_Z(\M)$.  It
sits on the blowup $B_Z(X)$, and is simply the preimage $b^*\M$ modulo
its submodule of sections supported on $D$, where $b$ is the blowup map
and $D$ is the exceptional divisor.  (In \cite{\KT, (2.1), \p.181},
this transform was denoted by $b^\dagger\M$.)  The transform $P_Z(\M)$
is simply the ``specialization'' of $\M$ to the completed normal cone
$P_Z$, defined via a version of the specialization process of
\cite{\Fulton, Ch.~5}: let $\wh X$ be the relative affine line $\Spec
\O_X[t]$, and $\wh \M$ the pullback of $\M$; embed $X$ in $\wh X$ as
the zero section, form the blowup $B_Z(\wh X)$ of $\wh X$ along $Z$,
and set
	$$P_Z(\M):=B_Z\big(\wh\M\,\big)_Z$$
 where the right hand side is the restriction to the preimage of $Z$;
this preimage is the exceptional divisor of $B_Z(\wh X)$, and is equal
to $P_Z$.  The transform $P^{1,2}_Z(\M)$ is defined similarly.  In
place of $\wh B_Z$, we use the joint blowup $\wh B_{1,2}=B_{1,2}(\wh
X)$ of $\wh X$ along two closed subschemes $Z_1$ and $Z_2$; by
definition, it is the scheme arising from the bigraded Rees algebra,
that is, the direct sum of the products of the powers of the ideals of
$Z_1$ and $Z_2$.  Although $\wh B_{1,2}$ is equal to the ordinary
blowup of $\wh X$, the center is the sum $Z_{1,2}$ of $Z_1$ and $Z_2$
in $\wh X$, not their sum $Z$ in $X$; nevertheless, $P^{1,2}_Z(\M)$ is
still defined via restriction to the preimage of $Z$.

The transforms $P_Z(\M)$ and $P^{1,2}_Z(\M)$ have dimension at most
$\dim\M$, and when $r$ is given, the formations of the cycles
$[P_Z(\M)]_r$ and $[P^{1,2}_Z(\M)]_r$ are additive in $\M$ if
$\dim\M\le r$.  Hence, given a cycle $\S$ on $X$, it is reasonable to
define transforms $P_Z(\S)$ and $P^{1,2}_Z(\S)$ as follows: if $\bold
S=[S]$ where $S$ is integral with dimension $r$, then the transforms
are the cycles $[P_Z(\O_S)]_r$ and $[P^{1,2}_Z(\O_S)]_r$; next, extend
this definition by linearity.  It follows that
  $$[P_Z(\M)]_r=P_Z([\M]_r)\and[P^{1,2}_Z(\M)]_r=P^{1,2}_Z([\M]_r).$$
 For details, see (4.1)--(4.2) and (5.1)--(5.2).

The preceding considerations suggest extending the definitions of the
various Buchsbaum--Rim multiplicities from a module $\M$ to an
$r$-cycle $\S$, and from certain intersection numbers to analogous
direct image cycles of arbitrary dimension, which we call ``twisted
Segre classes.''  To be more precise, let's fix an arbitrary invertible
sheaf $\L$; it plays the role here, played by $\O_X(1)$ in \cite{\KT}.
In (4.4), we define the $i$th {\it twisted Segre\/} class by the formula,
	$$s^i(Z,\L)(\S):=p_*c_1\L(1)^iP_Z(\S),$$
 where $p\:P_Z\to Z$ is the natural map and $\L(1)$ stands for
$p^*\L\ox \O(1)$.  For example, $s^i(Z,\O_X)[X]$ is the usual Segre
class of $Z$ in $X$.

By construction, $s^i(Z,\L)(\S)$ is a rational equivalence class on
$Z\cap S$ where $S$ is the support of $\S$.  In fact, $s^i(Z,\L)(\S)$
is, in a natural way, a sum of classes $s_W^i(Z,\L)(\S)$ supported on
the images $W$ of the components of the cycle $P_Z(\S)$.  Not
surprisingly, these $W$ play an important role, especially in the
proofs of the subtler properties of positivity in Section 7.  So, if a
subset of $Z\cap S$ is a $W$, then we say that it's {\it distinguished}
by $(Z,\S)$, or by $(Z,\M)$ if $\S=[\M]_r$.

Distinguished subsets are studied in Section 3.  None need exist; in
other words, the cycle $P_Z(\S)$ may vanish although $\S$ doesn't
vanish.  Moreover, a distinguished subset may exist although no
component of $Z\cap S$ is distinguished.  However, if an irreducible
subset $W$ of $Z\cap S$ satisfies the chain condition,
	$$\dim W+\cod(W,S)=\dim S,$$
 then $W$ is contained in a set $W_1$ distinguished by $(Z,\S)$.
Moreover, if $\S$ is prime, then $W$ itself is distinguished if and
only if both this chain condition obtains and the ideal
  $\I_Z$ has {\it maximal analytic spread} on $\O_S$ at the generic point
$w$ of $W$.
  The latter condition means simply that the subset
$\{w\}$ of the localization $X_w:=\Spec\O_{X,w}$ is distinguished by
the localization $(Z_w,\O_{S,w})$.  This definition is clearly
equivalent to the original definition of Northcott and Rees \cite{\NR,
Dfn.~3, \p.199}.

  If $\S=[\M]_r$, then we also say that $\I_Z$ has maximal analytic
spread on $\M$ at $w$.  If also $X$ is a local scheme with $w$ as its
closed point (so $W=\{w\}$), then this condition of maximal analytic
spread is detected by a polynomial, in $n$ say, of degree at most $r$.
The coefficient of $n^r/r!$ is equal to the intersection number,
	$$e(\I_Z,\M):=\int c_1\O_{P_Z}(1)^{r}[P_Z(\M)]^w,$$
 where $[P_Z(\M)]^w$ is the part of the fundamental cycle lying over
$w$.  It follows that $e(\I_Z,\M)$ is nonzero if and only if $\I_Z$ has
maximal analytic spread on $\M$ at $w$.  Moreover, $e(\I_Z,\O_X)$ is just the
multiplicity treated by Achilles and Manaresi \cite{\AM}, and if
$Z\cap\Supp\M=\{w\}$, then $e(\I_Z,\M)$ is just the usual Samuel
multiplicity.

One of our main theorems is the additivity theorem (4.6).  It concerns
a closed subscheme $W$ of $X$ containing $Z$, and the residual scheme
$R$ of the exceptional divisor $D$ in the preimage $b^{-1}W$ on the
blowup $B_Z(X)$.  The theorem asserts the following formula relating
rational equivalence classes on $W$:
	$$s^i(W,\L)(\S)=s^i(Z,\L)(\S)+b_*s^i(R,\L(1))B_Z(\S).$$
 This formula generalizes the additivity formula
for Buchsbaum--Rim multiplicities \cite{\KT, (6.7b)(i), \p.205}.  When
$\L=\O_X$, the formula recovers the residual intersection formula, and
gives an interpretation of its lower degree terms in terms of twisted
Segre classes.

The $(i,j)$-th \dfn{mixed twisted Segre\/} class is defined in (5.4) by
the formula,
	$$s^{i,j}(Z_1,\L_1;Z_2,\L_2)(\S)
	:=p^{1,2}_*c_1\L_1(1,0)^ic_1\L_2(0,1)^jP^{1,2}_Z(\S),$$
 where $p^{1,2}\:P^{1,2}_Z\to Z$ is the restriction of the joint blowup
map $\hat b\:\wh B_{1,2}\to\wh X$, where $\L_1$ and $\L_2$ are two
given invertible sheaves on $X$, and where $\L_1(1,0)$ and $\L_2(0,1)$
stand for their pullbacks twisted by the two tautological sheaves.

Another main theorem (5.5) asserts the following mixed-class formula,
which generalizes the mixed-multiplicity formula announced at the end
of \cite{\KT, (9.10)(ii), \p.225}: for any $n$,
 $$s^{n}(Z_1+Z_2,\L_1\otimes\L_2)(\S)=\sum_{i+j=n}\binom ni
	    s^{i,j}(Z_1,\L_1;Z_2,\L_2)(\S).$$
 The right side is obviously equal to
	$p^{1,2}_*c_1(\L_1\otimes\L_2(1,1))^nP^{1,2}_Z(\S)$.
 However, the latter is not obviously equal to the left side, because
$\wh B_{1,2}$ is the blowup along $Z_{1,2}$ not $Z$, and so a proof is
required.

The expansion formula (5.6.1) expresses the mixed classes in terms of
the pure classes.  Given it, the mixed-class formula is equivalent to
the special case of the additivity formula in which $Z$ is a divisor;
see (5.7).  However, by (4.7), this special case is equivalent to the
general case.  Thus, somewhat surprisingly, the mixed-class formula and
the additivity formula are essentially equivalent.

Section 6 considers suitably general `pseudo-divisors' $K$, and
determines the effect of replacing $\S$ by $K\cdot\S$ in both the pure
and the mixed twisted Segre classes.

Section 7 contains the deepest results, which concern positivity.  By
paying attention to the distinguished components, we improve some of
the statements in \cite{\KT}.  Moreover, we work in somewhat greater
generality: instead of assuming that $X$ is projective, $\L$ is
$\O_X(1)$, and so forth, we work with {\it nonnegative\/} and {\it
positive\/} rational equivalence classes {\bf s} and invertible sheaves
$\K$; we call {\bf s} \dfn{nonnegative\/} (resp., \dfn{positive\/}) and
write $\bold s \scq 0$ (resp., $\bold s\succ0$) if some multiple
$n\bold s$ with $n>0$ is represented by a nonnegative cycle (resp., by
a positive cycle).  Beware: a positive $\bold s$ can vanish in general,
but not when $X$ is projective if $\bold s$ is a class on the closed
fiber.  We call $\K$ \dfn{nonnegative\/} (resp., \dfn{positive\/}) and
write $\K \scq 0$ (resp., $\K\succ0$) if $c_1(\K)$ carries nonnegative
classes (resp., positive classes with positive dimension) into
nonnegative classes (resp., positive classes).  For example, $\K\scq0$
if $\K$ is generated by its global sections.  The main nonnegativity
result is (7.5).  It concerns a nonnegative cycle $\S$, and asserts,
for example, that $s^i(Z,\L)(\S)\scq0$ for $i\ge0$ if $\L|Z\scq0$ and
$\L(1)|D\scq0$; moreover, $s^{i,j}(Z_1,\L_1;Z_2,\L_2)(\S)\scq0$ for
$i,j\ge0$ if certain, more involved, conditions on $\L_1$ and $\L_2$
are satisfied.  Finally, the main positivity result (7.7) is one of the
main theorems.  It gives conditions, including the existence of certain
distinguished subsets, which imply that
$s^{i,j}(Z_1,\L_1;Z_2,\L_2)(\S)\succ0$ for $i$ and $j$ in a suitable
range.

\sectionhead 2 Transforms

\stp1
 Fix a scheme $X$, a closed subscheme $Z$, and a quasi-coherent
$\O_X$-module $\M$.  In this section alone, $X$ may be arbitrary, and
$\M$ need only be quasi-coherent; starting with the next section, we'll
assume, usually implicitly, that $X$ is Noetherian and that $\M$ is
coherent.  Form the blowup of $X$ along $Z$, and denote it by $B_Z(X)$
or $B_Z$ or $B$.  Denote the blowup map by $b\:B\to X$, the exceptional
divisor by $D$, and the tautological sheaf by $\O(1)$.

Set $U:=X-Z$, and denote the canonical embedding, which identifies $U$
with $B-D$, by $i\:U\to B$.  For short, denote the ideal $\I_Z$ of $Z$
in $X$ by $\I$.  Finally, let $W$ be an arbitrary closed subscheme of
$X$, and for short, denote its ideal  $\I_W$ by $\J$.

\art 2 Two transforms

Define a transform of $\M$ on $B_Z$  by
	$$B_Z(\M):=\Im(b^*\M\to i_*\M_U).$$
 This transform will play an important role, and we'll call it the {\it
proper transform}.  It sits in a short exact sequence,
      $$0\to B_Z(\M)(1)@>\delta>> B_Z(\M)\to B_Z(\M)_Z\to 0.$$
 Indeed, the section of $\O_B(-1)$ defining $D$ induces an isomorphism
from $(i_*\M_U)(1)$ to $i_*\M_U$.  So its restriction $\delta$ is an
injection.  The cokernel of $\delta$ is the restriction $B_Z(\M)_Z$
over $Z$ because $D=b^{-1}Z$.

Define another transform of $\M$ on $B$  by
	$$R_{Z,W}(\M):=\Im\big(B_Z(\M)_W\to B_Z(\M)_W(-1)\big).$$
 This transform will serve as a temporary device.  Its twist sits in the
short exact sequence,
 $$0\to R_{Z,W}(\M)(1)\to B_Z(\M)_W \to B_Z(\M)_{Z\cap W}\to0,
	\tgs2.1 $$
 which arises from the one above via restriction over $W$.

Obviously, $B_Z(\O_X)_W$ is equal to the structure sheaf of $b^{-1}W$;
hence, $R_{Z,W}(\O_X)$ is the structure sheaf of a closed subscheme
$R$, and \Cs2.1) yields this short exact sequence,
 $$0\to \O_R(1)\to \O_{b^{-1}W}\to \O_{D\cap b^{-1}W}
\to 0.\tgs2.2
 $$
 When $W$ contains $Z$, then $R$ is the `residual' scheme of $D$ in
$b^{-1}W$, which is defined by the equation $b^{-1}W=D+R$ (where the sum
of two schemes is defined by the product of their ideals).

On the other hand, $B_Z(\O_W)$ is the structure sheaf of the proper
transform of $W$, which is equal to the blowup of $W$ along $Z\cap W$.
Correspondingly, $B_Z(\M_W)$ is equal to $B_{Z\cap W}(\M_W)$, where
$Z\cap W$ and $W$ and $\M_W$ play the normal roles of $Z$ and $X$ and
$\M$.  Moreover, clearly, there are natural surjections,
 	$$B_Z(\M)_W\onto R_{Z,W}(\M)\onto B_Z(\M_W).$$

\art 3 The corresponding graded modules

Since the blowup $B$ is the `Proj' of the Rees algebra $\bigoplus\I^p$,
the sheaves and maps on $B$ in \Cs2) must arise from natural graded
modules and maps.  And indeed, the map from $b^*\M$ to $i_*\M_U$ arises
from the natural map from $\bigoplus\I^p\ox\M$ to $\bigoplus\M$, as can
be checked easily on principal open sets; hence, $B_Z(\M)$ arises from
$\bigoplus\I^p\M$, and the injection in \Cs2.1) arises from the
inclusion of $\I^{p+1}\M$ into $\I^p\M$.

The preimage $b^{-1}W$ is defined by the graded ideal $\bigoplus
\J\I^p$.  So $B_Z(\M)_W$ arises from the graded module,
 $$\bigoplus \I^p\M/\J\I^p\M. $$
 It follows that $R_{Z,W}(\M)$ arises from the graded module whose $p$th
piece is
	 $$\I^p\M/(\I^p\M\cap \J\I^{p-1}\M).\tgs3.1 $$
 Finally, $\M_W:=\M/\J\M$; hence, $B_Z(\M_W)$ arises from the graded
module whose $p$th graded piece is
	      $$\I^p\M/(\I^p\M\cap\J\M).\tgs3.2$$
 Moreover, the surjection from $R_{Z,W}(\M)$ onto $B_Z(\M_W)$ arises from
the graded map whose $p$th component is the obvious surjection from
\Cs3.1) onto \Cs3.2).

\lem4
 The natural surjection is an isomorphism,
 	 $$R_{Z,W}(\M)\risom B_Z(\M_W),\tgs4.1$$
 if, locally on $X$, there is a $p_0$ such that, for all $p\ge p_0$,

	$$\I^p\M\cap\J\M\subseteq\I^{p-1}\J\M.\tgs4.2$$
  The converse holds if $X$ is locally Noetherian and $\M$ is coherent.
  \pf
  A graded module ${\N}$ over the Rees algebra $\bigoplus\I^p$ gives rise
to the null sheaf on $B$ if ${\Cal{N}}_p=0$ for $p\gg0$ locally on $X$.
The converse holds if $\I$ and ${\N}$ are locally finitely generated.
Take ${\N}$ to be the kernel of the surjection from \Cs3.1) onto
\Cs3.2), and the assertion results.
 \enddemo

\art 5
  A third transform

Let $\wh X$ be the relative affine line $\Spec \O_X[t]$, and $\wh \M$
the pullback $\M[t]$.  View $X$ as embedded in $\wh X$ as the principal
divisor $\{t=0\}$.  Form the blowup  $B_Z(\wh X)$, and denote it as well
by $\wh B_Z$ or $\wh B$.  Note that the proper transform of $X$ in
$\wh B$ is equal to the blowup $B_Z(X)$ and that $\wh
\M_W=\M_W$; it follows that
	$$B_Z(\wh\M_W)= B_Z(\M_W).\tgs5.1$$

Denote the preimage of $W$ in $\wh B_Z$ by $P_W$.  For example, $P_Z$
is the exceptional divisor.  Define a transform of $\M$ on $P_W$ by
	$$P_W(\M):= B_Z\big(\wh\M\big)_W.$$
 This transform will play a leading role.  It arises from the
graded $\O_X[t]$-module with trivial $t$-action, whose $p$th piece
is the direct sum,
	$$\multline
	     \M/(\I+\J)\M\oplus\I\M/(\I+\J)\I\M\oplus\cdots\\
	\oplus\I^{p-1}\M/(\I+\J)\I^{p-1}\M\oplus\I^p\M/\J\I^p\M.
	\endmultline$$
 Indeed, this description is straightforward to check (if tedious).
Notice that $P_Z$ is the completed normal cone of $Z$ in $X$; also, if
$W\subseteq Z$, then $\I+\J=\J$, and so $P_W$ is the restriction over
$W$ of $P_Z$.

\lem6
 Keep the setup of\/ \Cs5).  Then
$R_{Z,X}\big(\wh\M\big)=B_Z(\M)$, and there is a short exact sequence,
	 $$0\to B_Z(\M)(1)\to P_X(\M)\to P_Z(\M)\to0.$$
 \pf
 Apply  \Cs4) with $\wh X$, $Z$, $X$ and
$\wh\M$ for $X$, $Z$, $W$ and $\M$.  The ideal of $X$ in $\wh
X$ is $(t)$; so, that of $Z$ is $(\I,t)$.  Hence, the left hand
side of \Cs4.2) becomes
	$$(\I,t)^p\M[t]\cap t\M[t].  $$
 Clearly, this intersection is equal to $(\I,t)^{p-1}t\M[t]$.  Therefore,
\Cs4) and \Cs5.1) yield the first assertion.  The second
assertion now follows from \Cs2.1).  \enddemo

\art 7 Joint blowups

Given two closed subschemes $Z_1$ and $Z_2$ of $X$, the `joint' blowup
$B_{1,2}$ or $B_{1,2}(X)$ is, by definition, the scheme arising from the
bigraded algebra $\bigoplus\I_{1}^p\I_{2}^q$ where $\I_k$ is the ideal
of $Z_k$.  It is not hard to identify $B_{1,2}$ as a repeated blowup:
first form the blowup $B_2$ of $X$ along $Z_2$; then form the blowup of
$B_2$ along the preimage of $Z_1$; that second blowup is equal to
$B_{1,2}$.  Of course, $B_{1,2}$ is, similarly, equal to the blowup of
$B_1$ along the preimage of $Z_2$.

On $B_{1,2}$, there is a `tautological' sheaf $\O(1,1)$; it arises from
the bihomogeneous ideal $\bigoplus\I_1^{p+1}\I_2^{q+1}$, and is equal to the
tensor product of the tautological sheaves on the $B_k$.  On $B_{1,2}$,
the preimage of $Z_k$ is a divisor $D_k$.  Set $D_{1,2}:=D_1+D_2$ and
take $Z:=Z_1+Z_2$; then $D_{1,2}$ is a divisor, it is equal to the
preimage of $Z$, and its ideal is equal to $\O(1,1)$.  Hence, there is
a natural map from $B_{1,2}$ to the blowup $B$ of $X$ along $Z$.  This
map has a natural inverse, because $B_{1,2}$ can be viewed as a
repeated blowup and because the pullback of each $Z_k$ to $B$ is a
divisor since the pullback of their sum $Z$ is a divisor.  Thus
$B_{1,2}$ and $D_{1,2}$ and $\O(1,1)$ are equal to $B$ and $D$ and
$\O(1)$.

In the current context, the transform $B_Z(\M)$ will be denoted by
$B_{1,2}(\M)$.  Clearly, it arises from the bigraded module
$\bigoplus\I_{1}^p\I_{2}^q\M$.  Moreover, $R_{Z,W}(\M)$ is equal to
$B_{1,2}(\M_W)$ if, locally on $X$, for $p,q\gg 0$,
 $$\I_1^p\I_2^q\M\cap \J\M\subseteq \I_1^{p-1}\I_2^{q-1}\J\M,
	\tgs7.1 $$
 and the converse holds if $X$ is locally Noetherian and $\M$ is
coherent.

The theory in \Cs5) generalizes naturally via the joint blowup $\wh
B_{1,2}$ of the relative affine line $\wh X$; beware, however, this
joint blowup is equal to the ordinary blowup along the sum $Z_{1,2}$ of
$Z_1$ and $Z_2$ in $\wh X$, not along their sum $Z$ in $X$.  Denote the
proper transform of $\wh\M$ on $\wh B_{1,2}$ by
$B_{1,2}\big(\wh\M\big)$.  Then the proper transform of $X$ in $\wh
B_{1,2}$ is equal to the blowup $B$ of $X$ along $Z$ because $Z_{1,2}$
meets $X$ in $Z$, and so, since $\wh \M_W=\M_W$, it follows that
  $$B_{1,2}\big(\wh\M_W\big)=B_{1,2}(\M_W).\tgs7.2$$
 Furthermore, set
	$$P^{1,2}_W(\M):= B_{1,2}\big(\wh\M\,\big)_W.$$

For each $p,q$, there is a natural filtration,
	$$\I_1^p\I_2^q=\I^{p,q}_{0}\subseteq\I^{p,q}_{1}\subseteq
	\cdots\subseteq \I^{p,q}_{p+q}=\cdots=\O_X.$$
 where $\I^{p,q}_\nu$ is the sum of all products $\I_1^i\I_2^j$ for
$i\le p$ and $j\le q$ and for $i+j\ge p+q-\nu$.  We'll now establish
the following inclusions for $\nu \ge 1$:
 $$\I_1\I_2\I^{p,q}_\nu \subseteq \I^{p,q}_{\nu -1}\hbox{ and }
 \I^{p,q}_{\nu}\subseteq \I^{p-1,q-1}_{\nu -1} \tgs7.3$$

 To prove these inclusions, it suffices to prove that their right hand
sides contain respectively $\I_1^{i+1}\I_2^{j+1}$ and $\I_1^i\I_2^j$
for $i\le p$ and $j\le q$ and $i+j\ge p+q-\nu$.  First assume
$i+j>p+q-\nu $.  Then $i+j\ge p+q-(\nu-1)$; so $\I_1^{i+1}\I_2^{j+1}$
is contained in the right side of the first inclusion in question.
Moreover, $(i-1)+(j-1)\ge (p-1)+(q-1)-(\nu -1)$; so $\I_1^i\I_2^j$ is
contained in that of the second.  Next assume $i<p$.  Then $(i+1)+j\ge
p+q-(\nu-1)$ and $i+1\le p$ and $j\le q$; so $\I_1^{i+1}\I_2^{j+1}$ is
contained in the right side of the first inclusion in question.
Moreover, $i\le p-1$ and $j-1\le q-1$ and $i+(j-1)\ge (p-1)+(q-1)-(\nu
-1)$; so $\I_1^i\I_2^j$ is contained in that of the second.  By
symmetry, the two inclusions hold when $j<q$.  Finally, the three cases
just considered are exhaustive (although not disjoint), because
$\nu\ge1$.

\lem8
 In the setup of\/ \Cs7), the transform $P_Z^{1,2}(\M)$ arises from the
bigraded $\O_X[t]$-module with trivial $t$-action, whose $(p,q)$-th
piece is the direct sum,
 $$(\I^{p,q}_0\M/\I_1^{p+1}\I_2^{q+1}\M)
 \oplus(\I^{p,q}_1\M/\I^{p,q}_0\M)
 \oplus\cdots\oplus(\I^{p,q}_{p+q}\M/\I^{p,q}_{p+q-1}\M).
\tgs8.1$$
 \pf
 The joint blowup $\wh B_{1,2}$ is defined using the ideals
$\?\I_k:=(\I_{Z_k},t)$, whereas  the ideal of $Z$ in $\wh X$ is
$\?\I:=(\I_1\I_2,t)$.  Obviously,
 $$\?\I_1^p\?\I_2^q=\I^{p,q}_0\oplus \I^{p,q}_1t\oplus
	\I^{p,q}_2t^2\oplus\cdots.\tgs8.2$$
 Consider the $\nu $th homogeneous component with respect to $t$ of
$\?\I\,\?\I_1^p\?\I_2^q$.  Clearly, for $\nu =0$, this component is
$\I_1^{p+1}\I_2^{q+1}$.  For $\nu>0$, the component is equal to the sum
$\I_1\I_2\I^{p,q}_\nu +\I^{p,q}_{\nu-1}$.  The latter is equal to
$\I^{p,q}_{\nu -1}$ by the first inclusion of \Cs7.3).  Hence,
 $$\?\I\,\?\I_1^p\?\I_2^q=\I_1^{p+1}\I_2^{q+1}\oplus\I^{p,q}_0t
	\oplus \I^{p,q}_1t^2 \oplus \cdots.\tgs8.3$$
 Clearly, the bigraded pieces of $\?\I_1^p\?\I_2^q\wh\M$ and
$\?\I\,\?\I_1^p\?\I_2^q\wh{\M}$ are obtained from \Cs8.2) and
\Cs8.3) simply by replacing $t^\nu$ by $\M t^\nu$.  Therefore, the
quotient,
 $$\?\I_1^p\?\I_2^q\wh{\M}/\?\I\,\?\I_1^p\?\I_2^q\wh{\M},$$
 is equal to \Cs8.1).  Hence, the assertion follows from the discussion
in \Cs7).\enddemo

\lem9
 In the setup of\/ \Cs7), there is a short exact sequence on $\wh
B_{1,2}$:
 $$0\to B_{1,2}(\M)(1)\to  P^{1,2}_X(\M)\to  P^{1,2}_Z(\M)\to 0.$$
 \pf
 The asserted exact sequence will arise from \Cs2.1) after we have
identified $B_{1,2}(\M)$ as $R_{Z_{1,2},X}(\wh\M)$.  To do so, it
suffices, because of \Cs7.2), to prove
	$$B_{1,2}\big(\wh\M_X\big)=R_{Z_{1,2},X}(\wh\M).$$
 To prove it, apply the criterion given in \Cs7) with $\wh X$,
$Z_{1,2}$, $X$ and $\wh\M$ for $X$, $Z$, $W$ and $\M$.  The ideal
of $X$ in $\wh X$ is $(t)$; so, that of $Z_k$ is $(\I_k,t)$.  Hence,
the left hand side of \Cs7.1) becomes
     $$\I_1^p\I_2^q\wh \M\cap t\wh \M.\tgs9.1$$
 The submodule $\I_1^p\I_2^q\wh \M$ of $\wh \M$ is,
obviously, the following:
 $$\I^{p,q}_0\M\oplus \I^{p,q}_1\M t\oplus\I^{p,q}_2\M t^2\oplus\cdots.
 $$
 So, by the second inclusion of \Cs7.3), the intersection \Cs9.1) lies
in $t\I_1^{p-1}\I_2^{q-1}\wh \M$ for $p,q\ge 1$.  Thus \Cs7.1)
holds, and the proof is complete.
 \enddemo

\sectionhead 3 Distinguished sets, maximal analytic spread

\stp1
 From now on, assume that $X$ is Noetherian with finite dimension, and
that the $\O_X$-module $\M$ is coherent.

For a moment, say that $X$ has dimension $r$.  Then the relative affine
line $\widehat X$ has dimension $r+1$, and its generic points belong to
the complement of $X$, which is viewed as the zero section.  Therefore,
the blowup $\widehat B$ of $\widehat X$ along $Z$ has dimension $r+1$,
and the preimage $P_W$ of $W$ in $\widehat B$ has dimension at most $r$,
where $W$ is a closed subscheme of $X$.  In particular, the exceptional
divisor $P_Z$ has dimension at most $r$.  It has dimension exactly $r$,
moreover, if and only if $Z$ has a (closed) point $z$ of codimension $r$
in $X$ (that is, $\dim\O_{X,z}=r$) by \cite{\KT, Lemma (3.2)(iii)} since
the dimension of $\widehat X$ at $z$ is 1 more than that of $X$;
furthermore, if so, then $z$ lies in the image of a component of $P_Z$
with dimension $r$.  Such an image will be said to be distinguished by
$(Z,X)$.

In general, a subset $W$ of $X$ will be said to be {\it
distinguished\/} by the pair $(Z,\M)$ if $W$ is the image of a
component $C$ of the support of $P_Z(\M)$ such that $\dim C=\dim\M$.
Notice that, in any case, $\dim C\le\dim\M$ by the reasoning above
applied after $X$ is replaced by the support $\Supp\M$ equipped with
the scheme structure defined by the annihilator.  Notice also that $W$
is closed, irreducible, and contained in $Z\cap\Supp\M$.  Furthermore,
$W$ will be said to be {\it distinguished\/} by a pair $(Z,S)$ where
$S$ is a given closed subscheme of $X$  if $W$ is distinguished by the pair
$(Z,\O_S)$.  Finally, the ideal $\I$ of $Z$ will be said to have {\it
maximal analytic spread\/} on $\M$ at a point $w$ if the subset $\{w\}$
of the localization $X_w:=\Spec\O_{X,w}$ is distinguished by the
localization $(Z_w,\M_w)$.

\prop2
 Let  $S$ be a closed subscheme of $X$.  Let
$W$ be an irreducible subset of $Z\cap S$, and $w$ its generic
point.  Then $W$ is distinguished by $(Z,S)$ if and only if {\rm (i)}
the ideal $\I$ of $Z$ has maximal analytic spread on $\O_S$ at $w$ and
{\rm (ii)} this equation holds:
	$$\dim W+\cod(W,S)=\dim S.\tgs2.1$$
 On the other hand, whenever \Cs2.1) holds, then $W$ is contained in a
set $W_1$ distinguished by $(Z,S)$ such that
	$$\dim W+\cod(W,W_1)=\dim W_1;\tgs2.2$$
 moreover, then $W$ is contained in a set $W_2$ distinguished by
$(Z_1,W_1)$ for any closed subscheme $Z_1$ of $X$ containing $W$.
 \pf
 Let $Z_w$, $S_w$ and so forth denote the localizations.  Clearly,
$P_Z(\O_S)$ is structure sheaf of a closed subscheme $P$ of $\wh B$,
and the formation of $P$ commutes with localization.  So $P_w$ may be
viewed as arising from $Z_w$ and $S_w$.  Hence, the reasoning in \Cs1)
yields $\dim P_w=\dim S_w$.

Let $C$ be an irreducible subset of $P$, and assume $C$ maps onto $W$.
Then,
	$$\dim C - \dim W = \dim C_w \le \dim P_w
	= \dim S_w \le \dim S - \dim W \tgs2.3$$
 Indeed, the first equality holds by \cite{\KT, Lemma (3.2)(ii)}; the
second equality was proved above; and the two inequalities are obvious.

First, assume $W$ is distinguished by $(Z,S)$.  Then, by definition,
there exists a $C$ such that $\dim C=\dim S$.  Hence the two
inequalities in \Cs2.3) are now equalities.  Thus $\dim C_w =\dim S_w$
and $\dim S_w = \dim S - \dim W$.  Hence (i) and (ii) hold.
Conversely, assume (i).  Then $P_w$ has a component $C'$ which maps
into $\{w\}$ and is such that $\dim C'=\dim S_w$.  Take the closure of
$C'$ in $P$ as $C$.  Then $C_w=C'$.  Hence, if \Cs2.1) holds also, then
\Cs2.3) yields $\dim C =\dim S$.  Therefore, $W$ is distinguished by
$(Z,S)$.

On the other hand, let $C_1$ be an irreducible subset of $P$, let $W_1$
be its image and $w_1$ the generic point of $W_1$, and assume $W_1\supseteq W$.
Then
	$$\dim C_1 - \dim W_1 = \dim C_{1,w_1} \and
	 \dim C_{1,w} - \dim W_{1,w} = \dim C_{1,w_1} \tgs2.4$$
 by \cite{\KT, Lemma (3.2)(ii)} applied to the map $C_1\to W_1$ and
again to the localization $C_{1,w}\to W_{1,w}$.

Recall that $\dim P_w=\dim S_w$.  Hence $P_w$ has a component $C'_1$ (not
necessarily contained in the closed fiber) such that $\dim C'_1=\dim
S_w$.  Take the closure of $C'_1$ in $P$ as $C_1$.  Then $C_{1,w}=C'_1$.
Hence,  \Cs2.4) yields
	$$\dim S_w  - \dim W_{1,w} = \dim C_1 - \dim W_1. \tgs2.5$$
 Now, $\dim S_w =\cod(W,S)$ and $\dim W_{1,w} = \cod(W,W_1)$ by
definition.  Moreover, $C_1\subseteq P$, and $\dim P\le\dim S$ by the
reasoning in \Cs1).  Hence
   $$\cod(W,S)-\cod(W,W_1) \le \dim S-\dim W_1.\tgs2.6$$

Assume \Cs2.1) also.  Then  \Cs2.6)  yields
	$$\dim W+\cod(W,W_1)\ge\dim W_1.$$
 However, the opposite inequality holds trivially.  Hence, \Cs2.2)
holds.  Moreover, \Cs2.6) is also an equality.  So \Cs2.5) now yields
$\dim C_1 = \dim S$.  Hence, $W_1$ is distinguished by $(Z,S)$.
Finally, the last assertion follows immediately from the one just
proved, but with $Z_1$ and $W_1$ for $Z$ and $S$.
 \enddemo

\art 3 Distinguished components

 Let $S$ be a closed subscheme with dimension $r$.  If $W$ is a
component of $Z\cap S$, then $W$ is distinguished by $(Z,S)$ if and only
if Equation \Cs2.1) holds.  Indeed, if it holds, then $W$ belongs to a
set $W_1$ distinguished by $(Z,S)$ by \Cs2); since $W_1$ is irreducible
and contained in $Z\cap S$, necessarily $W_1=W$.  The converse holds
directly by \Cs2).

If $W$ is a component of $Z\cap S$ with dimension $r$, then $W$ is a
component of $S$; it follows that $W$ is distinguished by $(Z,S)$.
Moreover, every component $W$ of $Z\cap S$ with dimension $r-1$ is
distinguished if $S$ is equidimensional with dimension $r$, because
Equation \Cs2.1) holds since the first summand is $r-1$ and the second
is at least 1.

Because the number of distinguished sets is finite, every component of
$Z\cap S$ is distinguished if every point  of $Z\cap S$ belongs to a
distinguished set.  The latter obtains if $S$ is biequidimensional (that
is, any two saturated chain of closed irreducible subsets have the same
length) because then Equation \Cs2.1) holds for every closed subset $W$
of $S$.

Every component of $Z\cap S$ is distinguished, obviously, if its closed
points are dense and if each one belongs to some distinguished subset.
The former condition obtains, by the Hilbert Nullstellensatz, if $Z\cap
S$ is of finite type over an Artin ring or over a Dedekind domain with
infinitely many primes.  The latter obtains, by \Cs2), if $S$ has
dimension $r$ at each closed point of $Z\cap S$.

\expl4
 Take as $X$ the spectrum of the polynomial ring in one variable $t$
over a discrete valuation ring with uniformizing parameter $\pi$.  Then
$X$ is regular of dimension 2.  However, the principal ideal $(\pi t-1)$
is maximal.  So it defines a divisor with dimension 0 and one point;
take it as $Z$.  Obviously, Equation \Cs2.1) fails for
$W:=Z$.  Hence, by \Cs2), there are no sets distinguished by $(Z,X)$; in
other words, $P_Z$ has dimension $1$, not $2$.

Take as $X$ the spectrum of a local domain that has dimension $3$ and is
not catenary.  Take as $Z$ an integral subscheme with dimension 1 and
codimension 1.  Then Equation \Cs2.1) fails for $W:=Z$; hence, $Z$
is not distinguished.  However, the equation holds for $W:=\{x\}$ where
$x$ the closed point of $X$.  Hence $\{x\}$ is the one and only
set distinguished by $(Z,X)$.

\prop5
 Assume that $X$ is a local scheme, and let $w$ be its closed point.
Assume that $Z$ is nonempty and $\M$ is nonzero, and set $r:=\dim\M$.
Let $Z_k$ be the $k$th infinitesimal neighborhood of $w$ in $Z$, and
let $\I$ be the ideal of $Z$.  Then, as a function of $n$, the length,
	$$\lambda(n):=\length\biggl(\biggl(\bigoplus_{i=0}^n
		\I^i\M/\I^{i+1}\M\biggr)\ox \O_{Z_k}\biggr),$$
 is eventually a polynomial of degree at most $r$, with equality if and
only if $\I$ has maximal analytic spread on $\M$ at $w$.  Moreover,
for $k\gg0$, the coefficient of
  $n^r/r!$
 is independent of $k$, and is equal to the intersection number,
 $$\int c_1\O_{\widehat B}(1)^r[P_Z(\M)]^w,\tgs5.1$$
 where $[P_Z(\M)]^w$ is the part of the fundamental cycle lying over
$w$.
 \pf
 It follows from the definition (2.5) that $P_{Z_k}(\M)$ is the
restriction over $Z_k$ of $P_Z(\M)$.  Therefore, the components of
$[P_Z(\M)]^w_r$ are equal to those of $[P_{Z_k}(\M)]_r$ for any $k$,
where the index $r$ indicates the part of dimension $r$; moreover, if
$k\gg0$, then they appear with the same multiplicities.  Such
components exist, obviously, if and only if $\I$ has maximal analytic
spread on $\M$ at $w$.

Since $P_{Z_k}(\M)$ is the above restriction, it arises from the graded
module whose $n$th piece is
	$$\biggl(\bigoplus_{i=0}^n
		\I^i\M/\I^{i+1}\M\biggr)\ox \O_{Z_k}.$$
 The latter has length $\lambda(n)$.  Therefore, by the familiar theory
of Hilbert polynomials (see \cite{\KT, Lemma (4.3)} for example),
$\lambda(n)$ is eventually a polynomial, whose degree is the dimension
of the support of $P_{Z_k}(\M)$.  Moreover, since this dimension is at
most $r$, the coefficient of $n^r/r!$ is equal to the asserted
intersection number \Cs5.1).
 \enddemo

\art 6 The generalized Samuel multiplicity

Let $e(\I,\M)$ denote the intersection number (3.5.1).  It is just the
multiplicity treated by Achilles and Manaresi; see \cite{\AM, (1.2) and
(1.3)}.  Moreover, $e(\I,\M)>0$ if and only if the fundamental cycle
$[P_Z(\M)]$ has $r$-dimensional components lying over $w$; in other words,
$e(\I,\M)>0$ if and only if  $\I$ has maximal analytic spread on $\M$ at $w$.

In \Cs5), $[P_Z(\M)]$ has dimension $r$ by the argument in \Cs1).
Hence, if $Z\cap\Supp\M=\{w\}$, then $\I$ has maximal analytic spread
on $\M$ at $w$.  If also $k$ is large enough so that $\M\ox\O_{Z_k}=
\M\ox\O_Z$, then $\lambda(n)$ is the ordinary Samuel polynomial of $\I$
on $\M$, and so $e(\I,\M)$ is the ordinary Samuel multiplicity.  Thus,
for an arbitrary $Z$ and $\M$, it is reasonable to call $e(\I,\M)$ the
\dfn{generalized Samuel multiplicity} of $\I$ on $\M$.  In the general
setup of \Cs1), it follows that the ideal $\I$ of $Z$ has maximal
analytic spread on $\M$ at $w$ if and only if $e(\I_w,\M_w)>0$.

\prop7
 Set $r:=\dim X$.  Then an irreducible subset $W$ of $Z$ is
distinguished by $(Z,X)$ if and only if either {\rm (i)} $\dim W=r$
(and $W$ is a component of both $Z$ and $X$) or {\rm (ii)} $\dim W<r$
and $W$ is the image of a component $\Gamma$ of $D$ with
$\dim\Gamma=r-1$.  In the former case, $W$ is equal to a component of
$P_Z$ disjoint from $D$; in the latter case, $W$ is the image of a
component $C$ of $P_Z$ such that $C\cap D=\Gamma$.
 \pf
 Set $P:=P_Z$.  Then $P =\Proj(\cG[u])$ where $\cG$ is the associated
graded algebra (or conormal algebra) of $Z$ in $X$ and where $\cG[u]$
is the polynomial algebra in one variable $u$ and is graded by total
degree; moreover, $D=\Proj(\cG)$ and the embedding $D\into P$ arises
from the canonical surjection $\cG[u]\onto\cG$.  The geometry here is
well known, and is the same, more generally, for any graded
$\O_Z$-algebra $\cG$ with $\cG_0=\O_Z$, with $\cG_1$ locally finitely
generated, and with $\cG$ generated by $\cG_1$.  Namely, the inclusion
$\cG\into \cG[u]$ gives rise to a central projection of $P$ onto $D$;
the center is the copy of $Z$ with ideal $\cG_1\!\cdot\,\cG[u]$, and the
blowup $V$ of $P$ along $Z$ is equal to the $\IP^1$-bundle
$\IP(\O_D\oplus \O_D(1))$.  So $P$ has two types of components: those
that are contained in $Z$, and those that aren't.  The former are also
components of $Z$, and are disjoint from $D$.  The latter correspond,
via proper transform, to components of $V$, which correspond, in turn,
via intersection (or via projection), to components of $D$.  Here,
proper transform preserves dimension by \cite{\KT, Lemma (3.2)(ii)},
and intersection decreases dimension by 1 because the operation is
inverted by forming the $\IP^1$-bundle.  The assertions follow directly.
 \enddemo

\sectionhead 4 Twisted Segre operators

\art 1 The transforms of a cycle

Let $\S$ be a cycle on $X$, and define natural analogs of the
transforms $B_Z(\M)$ and $P_W(\M)$ as follows: if $\S=[S]$ where $S$ is
integral with dimension $r$, then set
 $$\align
	B_Z(\S )&:=[B_Z(\O_S)]_r=[B_{Z\cap S}(S)]_r\\
	P_W(\S )&:=[P_W(\O_S)]_r=[B_{Z\cap S}(\wh S\,)_{W\cap S}]_r
 \endalign$$
 where the index $r$ indicates the part of dimension $r$; then extend
this definition by linearity.  Note that $B_{Z\cap S}(S)$ and $B_{Z\cap
S}(\wh S\,)_{W\cap S}$ have dimension at most $r$ by (3.1).  In
fact, if $S\subseteq Z$, then $B_{Z\cap S}(S)$ is empty; otherwise, it
is an integral scheme with dimension $r$.  On the other hand, $B_{Z\cap
S}(\wh S\,)_{W\cap S}$ may have several components with various
dimensions at most $r$.

Let $V$ be an arbitrary subset of $X$.  Say that $V$ is
\dfn{distinguished\/} by $(Z,\S)$ if $V$ is distinguished in the sense
of (3.1) by $(Z,S)$ where $S$ is the support of $\S$ given the induced
reduced structure.  Finally, given an arbitrary cycle $\T$ on a scheme
over $X$, denote the part (summand) of $\T$ formed by the components
whose generic points map {\it into} $V$ by ${\T}^V$ and the part formed
by the components that map {\it onto} $V$ by ${\T}\uprab V $.

\lem2
   A coherent $\O_X$-module $\M$ with dimension at most $r$ gives rise
to the following relation among cycles on $P_X$:
  $$[P_X(\M)]_r=[P_Z(\M)]_r+[B_Z(\M)]_r.\tgs2.1$$
 Moreover, the formation of each of the three cycles involved is
additive in $\M$, and those three cycles are equal, respectively, to
these:
 $$P_X([\M]_r)\and P_Z([\M]_r)\text{ and }B_Z([\M]_r).$$
 \pf
 Relation \Cs2.1) follows from the exact sequence of (2.6).  Now, by
construction, $B_Z(\M)$ has all its associated points in the complement
$U$ of the exceptional divisor, and $B_Z(\M)$ is equal to $\M$ on $U$.
It follows that $[B_Z(\M)]_r$ is equal to $B_Z([\M]_r)$.  Moreover,
therefore, the former transform is additive in $\M$ as the latter is.

By  construction, $B_Z(\wh\M\,)$ has all its associated points off the
(Cartier) divisors $P_Z$ and $P_X$.  Hence,
	$$[P_Z(\M)]_r=P_Z\cdot [B_Z(\wh \M)]_{r+1},
	\and [P_X(\M)]_r=P_X\cdot [B_Z(\wh \M)]_{r+1}.$$
 The asserted additivity in $\M$ follows since the two right hand sides
are additive by the first part of the proof.  Moreover, therefore, the
left hand sides are equal, respectively, to $P_Z([\M]_r)$ and
$P_X([\M]_r)$ because the latter transforms are additive in $\M$, and
equality holds by definition when $\M$ is the structure sheaf of an
integral subscheme $S$ of dimension $r$.  Alternatively, because of
\Cs2.1), the additivity of $[P_Z(\M)]_r$ may be concluded from that of
the other two transforms.
 \enddemo

\lemx{Key relations}3
 A cycle $\S$ on $X$ gives rise to the following relations among cycles
on $P_X$:
 $$\gather
 	P_X(\S)=P_Z(\S)+B_Z(\S);\tgs3.1\\
  P_Z(\S)=P_X(\S)^Z \text{ and }B_Z(\S)=P_X(\S)^{X-Z}.\tgs3.2
 \endgather$$
 Moreover, for any closed subscheme $W$ of $X$, the cycle $\S$ gives
rise to the following relation modulo rational equivalence on $P_W$:
	$$\hat hP_Z(\S)\uprab W =(D\cdot B_Z(\S))\uprab W \tgs3.3$$
 where $\hat h:=c_1\O_{\wh B_Z}(1)$ and $D$ is the exceptional divisor
of $B_Z$.
 \pf
 By linearity, we may assume that $\S=[S]$ where $S$ is integral
with dimension $r$.  Then, replacing $X$ by $S$, we have to prove
 $$\gathered
 [P_X]_r=[P_Z]_r+[B_Z]_r,\ [P_Z]_r=[P_X]_r^Z,\ [B_Z]_r=[P_X]_r^{X-Z},\\
	\and\hat h([P_Z]_r\uprab W )=(D\cdot[B_Z]_r)\uprab W .
 \endgathered\tgs3.4$$
 The first of these relations holds because of (2.6); in fact, (2.6)
yields the corresponding relation among divisors on $\wh B_Z$:
	$$P_X=P_Z+B_Z.\tgs3.5$$
 The next two relations follow because $B_Z$ has all its generic points
outside the exceptional divisor.

To prove the last relation, note that $\hat h=c_1\O(-P_Z)$.  Let $\T$
be a cycle on $\wh B_Z$, and work modulo rational equivalence on its
support.  Then $\hat h\T=-P_Z\cdot\T$.  Moreover, $P_X$ is principal;
so $P_X\cdot\T=0$.  Hence, \Cs3.5) implies that $\hat h\T=B_Z\cdot\T$.
Finally, take $\T:=[P_Z]_r\uprab W $ and note that $B_Z\cdot[P_Z]_r\uprab W
=[D]_r\uprab W $ because of (3.7).
 \enddemo

\art 4 Twisted Segre operators

 Fix an invertible sheaf $\L$ on $X$, and set
	$$\ell:=c_1(\L),\ h:=c_1\O_{B_Z}(1)
	\text{ and } \hat h:=c_1\O_{\wh B_Z}(1).$$
 Define the $i$th \dfn{twisted Segre operator\/} to be the following map
of degree $-i$ from cycles to rational equivalence  classes:
 $$\gathered
 s^i(Z,\L)\:\Cyc_r(X)\to \A_{r-i}(Z)\\
 s^i(Z,\L)(\S):=p_*(\ell+\hat h)^iP_Z(\S)
 \endgathered\tgs4.1
 $$
 where $r$ is arbitrary and $p\:P_Z\to Z$ is the restriction of $\hat
b\:\wh B_Z\to\wh X$.  Define the {\it total twisted Segre operator}
$s(Z,\L)$ to be the sum of the various $s^i(Z,\L)$.  So
	$$s(Z,\L)(\S)=p_*s(\L(1))P_Z(\S)\tgs4.2$$
 where $s(\L(1))$ is the usual total Segre operator $1/(1-(\ell+\hat h))$.
Set
	$$s(Z):=s(Z,\Cal O_X).$$
  Obviously, $s(Z)[X]$ is the usual Segre class of $Z$ in $X$.

By construction, $s(Z,\L)(\S)$ has support in $Z\cap S$ where $S$ is
the support of $\S$.  In fact, each $s^i(Z,\L)(\S)$ is, in a natural
way, a sum of classes $s_W^i(Z,\L)(\S)$ supported on the various sets
$W$ distinguished by $(Z,\S)$; this decomposition arises from that of
$P_Z(\S)$ into the pieces whose components map onto the $W$.

When $r$ is given, define for any coherent $\O_X$-module $\M$ such
that $r\ge\dim\M$,
	$$s(Z,\L)(\M):=s(Z,\L)([\M]_r)=p_*s(\L(1))[P_Z(\M)]_r$$
 where the second equation holds by \Cs2).  Obviously, $s(Z,\L)(\M)$
vanishes when $r>\dim\M$.

 Note  the following additivity formula:
 $$s^i(Z,\L\otimes \L_1)=\sum_{j}\binom{i}j s^j(\L)s^{i-j}(Z,\L_1).
 \tgs4.3
 $$
 It follows directly from the definition and the projection formula.

Note also that  $s^0(Z,\L)$ is given by the following formula:
	$$s^0(Z,\L)(\S)=\S^Z.\tgs4.4$$
 Indeed, by definition, $s^0(Z,\L)(\S)=p_*P_Z(\S)$.  Now,
$P_Z(\S)=P_X(\S)^Z$ by \Cs3.2).  Moreover, clearly, if $\wh\S$ denotes
(of course) the pullback of $\S$ to $\wh X$, then
	$$p_*P_X(\S)=p_*(P_X\cdot B_Z(\wh\S\,))=X\cdot p_*B_Z(\wh\S\,)
	=X\cdot\wh\S=\S,$$
 and the asserted formula  \Cs4.4) follows.

Standard intersection theory and standard properties of blowups yield,
in a straightforward fashion, the following two functorial properties of
$s(Z,\L)$ with respect to a map $q\:X'\to X$:

(a) (Proper pushforward).  If $q$ is proper, then we have the following
identity of operators from $\Cyc(X')$ to $\A(Z)$:
 $$q_*s(q^{-1}Z,q^*\L)(\S')=s(Z,\L)q_*(\S'). $$

(b) (Flat pullback).  If $q$ is flat and if, at every generic point of
every fiber, the residue field extension has the same transcendence
degree, then we have the following identity of operators from $\Cyc(X)$
to $\A(q^{-1}Z)$:
 $$q^*s(Z,\L)(\S)=s(q^{-1}Z,q^*\L)q^*(\S).$$

\art 5 Blowup formulas

In terms of the blowup map $b\:B_Z\to X$ and its exceptional divisor
$D$, the twisted Segre operators are given by the following {\it blowup
formulas}:
 $$\align
  s^i(Z,\L)(\S)&=\ell^i\S^Z+ b_*\frac{(\ell+h)^i-\ell^i}{h}D\cdot
		B_Z(\S) \tgs5.1\\
	&=\ell^i\S^{Z}+\sum_{j+k=i-1}b_*\ell^j(\ell+h)^kD\cdot B_Z(\S)
		\tgs5.2\\
	&=\ell^i\S^{Z}-b_*\bigl((\ell-D)^i-\ell^i\bigr)B_Z(\S);
		\tgs5.3\\
  s(Z,\L)(\S)&=s(\L)\S^Z+ b_*s(\L)s(\L(1))D\cdot B_Z(\S).\tgs5.4\\
 \endalign$$
 Here, the divisor $D$ is transversal to the $r$-cycle $B_Z(\S)$, and
their intersection is an $(r-1)$-cycle.  The operator
$(\ell-D)^i-\ell^i$ is to be interpreted as follows: expand it as a sum
of terms $\ell^jD^k$ where $k>0$, and evaluate $\ell^jD^k$ on
$B_Z(\S)$, getting a class in $\A(D\cap B_Z(S))$, by first forming $D\cdot
B_Z(\S)$ and then applying $\ell^jD^{k-1}$.

 Those formulas are valid modulo rational equivalence on $Z\cap S$ where
$S$ is the support of $\S$, and they are obviously all equivalent.  The
next two formulas are valid on $S$, and they follow from the others
because, clearly, $b_*B_Z(\S)=\S^{X-Z}$.
 $$\align
  s^i(Z,\L)(\S)&=\ell^i\S-b_*(\ell-D)^iB_Z(\S);\tgs5.5\\
  s(Z,\L)(\S)&=s(\L)\S-b_*s(\L(1))B_Z(\S).\tgs5.6
 \endalign$$

To prove \Cs5.1), rewrite \Cs4.1) using the identity,
	$$(\ell+\hat h)^i
	=\ell^i+\frac{(\ell+\hat h)^i-\ell^i}{\hat h}{\hat h}.$$
 Apply \Cs4.4) and the projection formula to the first summand, and
apply \Cs3.3) to the second.  The result is \Cs5.1).

Consider the special case where $Z$ is a divisor $C$ in $X$.  Then
$B_Z=X$ and $D=C$ and $B_Z(\S)=\S^{X-C}$.  So \Cs5.4) and \Cs5.6)
become
 $$\align
  s(C,\L)(\S)&=s(\L)(\S^C)+ s(\L)s(\L(-C))C\cdot \S^{X-C}\tgs5.7\\
	&=s(\L)\S-s(\L(-C))\S^{X-C},\tgs5.8
 \endalign$$
and they are valid on $C\cap S$ and $S$ respectively.

\thmx {Additivity}6
 Let $W$ be a closed subscheme of $X$ containing the subscheme $Z$, and
let $R$ be the residual scheme on $B_Z$ of $D:=b^{-1}Z$ in $b^{-1}W$.
Then we have the following relation among operators from $\Cyc(X)$ to
$\A(W)$:
   $$s(W,\L)(\S)=s(Z,\L)(\S)+b_*s(R,\L(1))B_Z(\S).\tgs6.1$$
 \pf
 First note that $R$ is equal to the residual scheme $\wh R$ on $\wh
B_Z:=B_Z(\wh X)$ of $P_Z$ in $P_W$.  Indeed, by definition, $P_W=P_Z+\wh
R$.  Intersecting that equation with $B_Z$ yields $b^{-1}W=D+(\wh R\cap
B_Z)$.  So it remains to show that $\wh R\subseteq B_Z$.  The latter
may be checked locally on $\wh B_Z$.  So say $P_Z:\zeta=0$, and consider
the ideals $\I_{\wh R}$ and $\I_{B_Z}$.  Then $\zeta\I_{\wh R}$ is the
ideal of $P_W$ by definition, and $\zeta\I_{B_Z}$ is that of $P_X$ by
(2.6).  Now, $P_W\subseteq P_X$; so $\zeta\I_{\wh
R}\supseteq\zeta\I_{B_Z}$.  Since $\zeta$ is regular, therefore $\I_{\wh
R}\supseteq\I_{B_Z}$, as required.

To prove \Cs6.1), we may assume, by linearity, that $\S=[S]$ where
$S$ is integral, and then replace $X$ by $S$.  Let $\tilde b\:\wt
B\to\wh B_Z$ be the blowup map with center $R$, and $\wt D_R$ the
exceptional divisor.  Let $\wt P_Z$ and $\wt P_W$ and $\wt P_X$ be
the preimages of $P_Z$ and $P_W$ and $P_X$.  Then $\wt P_Z$ is a
divisor because $P_Z$ is, and $\wt P_X$ is principal because $P_X$ is
(because $X$ is).  By the projection formula for divisors, $\tilde
b_*[\wt P_X]_r^W=[P_X]_r^W$; hence the first equation of \Cs3.4)
yields
     $$\tilde b_*[\wt P_X]_r^W=[B_Z]_r^R+[P_Z]_r.\tgs6.2$$

By the first paragraph, $P_W=P_Z+R$.  So $\wt P_W=\wt P_Z+\wt D_R$.
Hence $\wt P_W$ is a divisor on $\wt B$.  Therefore, there is a (unique)
map $q\:\wt B\to\wh B_W$ such that $\wt P_W=q^{-1}\wh D_W$ where $\wh
B_W$ is the blowup of $\wh X$ along $W$ and where $\wh D_W$ is the
exceptional divisor.  Correspondingly, we have the following relation on
$\wt B$ among (the pullbacks of) the tautological $c_1\!$'s:
	$$\hat h_W=\hat h_Z+\hat h_R.\tgs6.3$$
 Since $q$ is an isomorphism over the complement of $W$, the projection
formula and the first equation of \Cs3.4) yield, as above,
	$$q_*[\wt P_X]_r^W=[\wh D_W]_r.\tgs6.4$$

Let $B^*$ be the proper transform of $B_Z$ under $\tilde b$.  So $B^*$
is equal to the blowup of $B_Z$ along $R$, and $\wt D_R\cap B^*$, which
is simply the preimage, $D_R$ say, of $R$ in $B^*$, is equal to the
exceptional divisor.  Since $\tilde b$ is an isomorphism off $R$,
obviously,
 $$[\wt P_X]_r=[B^*]_r+[\wt P_X]_r^W.\tgs6.5$$
 Since $\wt P_X$ is principal, $\wt P_X\cdot [\wt D_R]_r$ vanishes in
$\A(\wt D_R)$.  Hence, by the commutativity of intersection product,
$\wt D_R\cdot[\wt P_X]_r$ vanishes in $\A(\wt D_R)$.  So \Cs6.5) yields
	$$\wt D_R\cdot[B^*]_r=-\wt D_R\cdot [\wt P_X]_r^W$$
 in $\A(\wt D_R)$, so on $\A(\wt P_W)$.  Therefore, in $\A(\wt P_W)$,
 	$$\hat h_R[\wt P_X]_r^W = D_R\cdot [B^*]_r.\tgs6.6 $$

On $\wh B_Z$, set
	$\ell_Z:=\ell+\hat h_Z=c_1(\L(1)).$
 On $\wt B$,  \Cs6.3) yields the identity,
	$$(\ell+\hat h_W)^i=\ell_Z^{\,i}
	+\bigr((\ell_Z+\hat h_R)^i-\ell_Z^{\,i}\bigl).$$
 Apply the operators on both sides to $[\wt P_X]_r^W$, and push the
result forward into $\A(W)$.  On the left, push via $\A(\wh D_W)$; then
\Cs6.4) and the projection formula yield $s^i(W,\L)([X])$, which is,
after summing over $i$, the left side of \Cs6.1) because $\S=[X]$.

  On the right, push via $\A(P_W)$.  In $\A(P_W)$, \Cs6.2) and the
projection formula and \Cs6.6) yield
	$$\ell_Z^{\,i}[P_Z]_r+\ell_Z^{\,i}[B_Z]_r^R
	+\tilde b_*\frac{(\ell_Z+\hat h_R)^i
	-\ell_Z^{\,i}}{h_R}D_R\cdot [B^*]_r.$$
 That sum is, by the blowup formula \Cs5.1), equal this sum:
	$$\ell_Z^{\,i}[P_Z]_r+s^i(R,\L(1))[B_Z]_r.$$
 Pushed into $\A(W)$, this sum becomes, after summing further over $i$,
the right side of \Cs6.1) because $\S=[X]$.  Thus \Cs6.1) is proved.
 \enddemo

\art 7 A divisorial center

 Consider the additivity formula \Cs6.1) when $Z$ is a divisor $C$ in
$X$.  Then $R$ is the closed subscheme of $X$ defined by the equation
$W=C+R$.  Moreover, $b=1_X$ and $\O(1)=\O(-C)$ and $B_Z(\S)=\S^{X-C}$.
So \Cs6.1) becomes
	$$s(C{+}R,\L)(\S)=s(C,\L)(\S)+s(R,\L(-C))(\S^{X-C}).
	\tgs7.1$$

Conversely, the special case \Cs7.1) of \Cs6.1) implies the general
case.  Indeed, $b^{-1}W$ contains the exceptional divisor $D:=b^{-1}Z$,
and $\O(1)=\O(-D)$; moreover, $R$ is defined by $b^{-1}W=D+R$.  So
\Cs7.1), with $B_Z$ and $D$ for $X$ and $C$, yields
	$$s(b^{-1}W,b^*\L)B_Z(\S)=s(D,b^*\L)B_Z(\S)+s(R,\L(1))B_Z(\S).
	\tgs7.2$$
 Now, by linearity, we may assume that $\S$ is prime.  Applying $b_*$
to \Cs7.2) yields \Cs6.1) because of \Cs4)(a) if the support $S$ of
$\S$ does not lie in $Z$.  Otherwise, both sides of \Cs6.1) are equal
to $s(\L)\S$ by \Cs5.4) because $B_Z(\S)=0$.  Thus \Cs7.1) and
\Cs6.1) are equivalent.

In (5.7), we'll derive \Cs7.1) from two other formulas.  Then we'll
have a second proof of \Cs6.1).

\sectionhead 5 Mixed Segre operators

 \art 1 The joint transforms of cycles

Fix two closed subschemes $Z_1$ and $Z_2$ of $X$, and assume that
$Z_1+Z_2=Z$.  As in (2.7), form the blowups $B_k:=B_{Z_k}(X)$ and the
joint blowup $B_{1,2}:=B_{1,2}(X)$.  Let $D_k$ and $D$ be the pullbacks
of $Z_k$ and $Z$ to $B_{1,2}$.  Set $h_k:=c_1\O_{B_k}(1)$.  Finally,
set $\wh B_k:=B_{Z_k}(\wh X)$ and $\wh B_{1,2}:=B_{1,2}(\wh X)$ where
$\wh X$ is the relative affine line; also set $\hat h_k:=c_1\O_{\wh
B_{k}}(1)$ and $\hat h:=c_1\O_{\wh B_{1,2}}(1)$.

Let $\S$ be a cycle on $X$, and define natural analogs of the
transforms $B_{1,2}(\M)$ and $P^{1,2}_W(\M)$ as follows: if $\bold
S=[S]$ where $S$ is integral with dimension $r$, then set
 $$\align
	B_{1,2}({\bold{S}})&:=[B_{1,2}(\O_S)]_r=[B_{1,2}(S)]_r,\\
	P^{1,2}_W({\bold{S}})&:=[P^{1,2}_W(\O_S)]_r=[B_{1,2}(\wh S\,)_W]_r
 \endalign$$
 where $B_{1,2}(S)$ is the joint blowup of $S$ along $Z_1\cap S$ and
$Z_2\cap S$; then extend this definition by linearity.  Note that
$B_{1,2}(S)$ and $B_{1,2}(\wh S\,)_W$ have dimension at most $r$ by the
reasoning in (3.1).

The next two lemmas are modeled on (4.2) and (4.3), and may be proved
similarly, except for the last assertion of \Cs3); so, only its proof
will be given.  This proof requires a modified approach, because $\wh
B_{1,2}$ is the blowup along the sum $Z_{1,2}$ of $Z_1$ and $Z_2$ in
$\wh X$, not their sum $Z$ in $X$; so $P^{1,2}_Z$ is not necessarily a
divisor.  For the same reason, although there are the two proofs of
additivity in $\M$ of $[P_Z(\M)]_r$ in (4.2), only the second, the
alternative proof, carries over.

\lem2
  A coherent $\O_X$-module $\M$ with dimension at most $r$ yields
the following relation among cycles on $P^{1,2}_X$:
  $$[P^{1,2}_X(\M)]_r=[P^{1,2}_Z(\M)]_r+[B_{1,2}(\M)]_r.\tgs2.1$$
 Moreover, the formation of each of the three cycles involved is
additive in $\M$, and those three cycles are equal, respectively, to these:
 $$P^{1,2}_X([\M]_r)\and P^{1,2}_Z([\M]_r)\text{ and }B_{1,2}([\M]_r).$$
 \endgroup

\lemx{Key relations}3
 A cycle $\S$ on $X$ yields the following relations among cycles
on $P^{1,2}_X$:
 $$\gather
 	P^{1,2}_X(\S)=P^{1,2}_Z(\S)+B_{1,2}(\S).\tgs3.1\\
  P^{1,2}_Z(\S)=P^{1,2}_X(\S)^Z \text{ and }
	B_{1,2}(\S)=P^{1,2}_X(\S)^{X-Z}.\tgs3.2
 \endgather$$
 Moreover, $\S$ yields the following relations modulo rational
equivalence on $P^{1,2}_Z$:
 $$\hat h_kP^{1,2}_Z(\S)=D_k\cdot B_{1,2}(\S)\and
	\hat hP^{1,2}_Z(\S)=D\cdot B_{1,2}(\S).\tgs3.3$$
 \pf As noted before \Cs2), we need only prove \Cs3.3).
 By linearity, we may assume that $\S=[S]$ where $S$ is integral
with dimension $r$.  Then, replacing $X$ by $S$, we have to prove
 $$\hat h_k[P^{1,2}_Z]_r=D_k\cdot [B_{1,2}]\and
	\hat h[P^{1,2}_Z]_r=D\cdot [B_{1,2}]\tgs3.4$$
 Let $\wh D_k$ be the pullback to $\wh B_{1,2}$ of the exceptional
divisor of $\wh B_k$.  Then the pullback of $\hat h_k$ is equal to
$c_1\O(-\wh D_k)$.  So
	$$\hat h_k[P^{1,2}_Z]_r=-\wh D_k\cdot[P^{1,2}_Z]_r.$$
 Now, $P^{1,2}_X$ is principal.  Hence $P^{1,2}_X\cdot [\wh D_k]_r$
vanishes in $\A(\wh D_k)$.  Therefore, $\wh D_k\cdot [P^{1,2}_X]_r$
vanishes in $\A(\wh D_k)$, so in $\A(P^{1,2}_Z)$ too.  Consequently,
\Cs3.1) yields
	$$\hat h_k[P^{1,2}_Z]_r=\wh D_k\cdot [B_{1,2}].$$
 However, $\wh D_k\cap B_{1,2}=D_k$.  Hence, the first relation of
\Cs3.4) holds.  The second relation follows because $\hat h_1+\hat
h_2=\hat h$ and $D_1+D_2=D$.

\art 4 Mixed twisted Segre operators

 Fix two invertible sheaves $\L_1$ and $\L_2$ on $X$, and set
$\ell_k:=c_1(\L_k)$.  Define the $(i,j)$-th \dfn{mixed twisted Segre
operator\/} to be the following map of degree $-(i+j)$ from cycles to
rational equivalence classes:
 $$\gathered
 s^{i,j}(Z_1,\L_1;Z_2,\L_2)\:\Cyc_r(X)\to \A_{r-i-j}(Z)\\
 s^{i,j}(Z_1,\L_1;Z_2,\L_2)(\S)=
	p^{1,2}_*(\ell_1+\hat h_1)^{i}(\ell_2+\hat h_2)^{j}P^{1,2}_Z(\S)
 \endgathered\tgs4.1
 $$
 where $r$ is arbitrary, $Z:=Z_1+Z_2$ is the sum in $X$, and
$p^{1,2}\:P^{1,2}_Z\to Z$ is the restriction of $\hat b\:\wh
B_{1,2}\to\wh X$.  Obviously, the class $s^{i,j}(Z_1,\L_1;Z_2,\L_2)(\S)$
has support in $Z\cap S$ where $S$ is the support of $\S$.
 Obviously, the operators are symmetric:
$$s^{i,j}(Z_1,\L_1;Z_2,\L_2)(\S)=s^{j,i}(Z_2,\L_2;Z_1,\L_1)(\S).\tgs4.2$$
 When $Z_1$ is empty, then $\wh B_{1,2}=B_{Z_2}(\wh X)$ and so
  $$s^{i,j}(\emptyset,\L_1;Z_2,\L_2)(\S)=s^i(\L_1)s^j(Z_2,\L_2)(\S).$$
  The following formula is similar to (4.4.4), and its proof is similar
too:
	$$s^{0,0}(Z_1,\L_1;Z_2,\L_2)(\S)=\S^Z.\tgs4.3$$

The following {\it blowup formula\/} is similar to (4.5.1):
  $$\split
 &s^{i,j}(Z_1,\L_1;Z_2,\L_2)(\S)=\ell_1^{\,i}\ell_2^{\,j}\S^Z\\
  &\qquad+b_{1,2*}\frac{(\ell_1+h_1)^i-\ell_1^i}{h_1}
		\ell_2^{\,j}D_1\cdot B_{1,2}(\S)
	+b_{1,2*}\frac{(\ell_2+h_2)^j-\ell_2^j}{h_2}
		\ell_1^{\,i}D_2\cdot B_{1,2}(\S) \\
  &\qquad\qquad-b_{1,2*}\frac{(\ell_1+h_1)^i-\ell_1^i}{h_1}\,
 \frac{(\ell_2+h_2)^j-\ell_2^j}{h_2}D_1\cdot D_2\cdot\rlap{$B_{1,2}(\S).$}
   \endsplit \tgs4.4$$
 The proof is similar too.  Set $\ell_k^{\,\prime}:=\ell_k+\hat h_k$,
and rewrite \Cs4.1) using the identity,
	$$\ell_1^{\,\prime i}\ell_2^{\,\prime j}=\ell_1^{\,i}\ell_2^{\,j}
	+(\ell_1^{\,\prime i}-\ell_1^i)\ell_2^{\,j}
	+(\ell_2^{\,\prime j}-\ell_2^j)\ell_1^{\,i}
	+(\ell_1^{\,\prime i}-\ell_1^i)(\ell_2^{\,\prime j}-\ell_2^j).$$
 Apply \Cs4.3) and the projection formula to the first summand,
apply \Cs3.3) to the remaining summands, and observe that
	$$\hat h_1D_2\cdot B_{1,2}(\S)=c_1\O(-D_1)D_2\cdot B_{1,2}(\S)
		=-D_1\cdot D_2\cdot B_{1,2}(\S).$$
 The result is \Cs4.4).

When $r$ is given, define for any coherent $\O_X$-module $\M$
such that $r\ge\dim\M$,
	$$\align
   s^{i,j}(Z_1,\L_1;Z_2,\L_2)(\M):&=s^{i,j}(Z_1,\L_1;Z_2,\L_2)([\M]_r)\\
 &=p^{1,2}_*(\ell_1+\hat h_1)^{i}(\ell_2+\hat h_2)^{j}[P^{1,2}_Z(\M)]_r
	\endalign$$
 where the second equation holds by \Cs2).  Obviously,
$s^{i,j}(Z_1,\L_1;Z_2,\L_2)(\M)$
vanishes when $r>\dim\M$.

\thmx{The Mixed Operator Formula}5
 For any $n$,
 $$s^{n}(Z_1+Z_2,\L_1\otimes\L_2)(\S)=\sum_{i+j=n}\binom ni
	    s^{i,j}(Z_1,\L_1;Z_2,\L_2)(\S). \tgs5.1
 $$
  \pf
 Set $\ell:=\ell_1+\ell_2$ and note that $\hat h:=\hat h_1+\hat h_2$.
Then, in \Cs5.1), the sum is, in view of \Cs4.1), obviously equal to
	$$p^{1,2}_*(\ell+\hat h)^nP^{1,2}_Z(\S). \tgs5.2$$
 Now, thanks to \Cs3.1), the proof of (4.5.1) shows that \Cs5.2) is
equal to
    $$\ell^n\S^{Z}+b_*\frac{(\ell+\hat h)^n
    -\ell^n}{\hat h}D_{1,2}\cdot B_{1,2}(\S).\tgs5.3$$
 However, $B_{1,2}(\S)=B_Z(\S)$.  Therefore, \Cs5.3) is equal, by
(4.5.1), to the left side of \Cs5.1), and the proof is complete.
 \enddemo

\propx{The expansion formula}6
 The mixed operators may be expanded in terms of the ordinary ones:
  $$\split
  s^{i,j}(Z_1,\L_1;Z_2,\L_2)&(\S)=s^i(\L_1)s^j(Z_2,\L_2)(\S)\\
  +&b_{2*}s^j(\L_2(1))s^i(b_2^{-1}Z_1,b_2^*\L_1)B_2(\S).
  \endsplit \tgs6.1$$
 \pf
 The formula makes sense.  Indeed, $s^j(Z_2,\L_2)(\S)$ is a class on
$Z_2\cap S$, so on $Z\cap S$, where $S$ is the support of $\S$.  Now,
$b_2\:B_2\to X$ is the blowup along $Z_2$, and
$\L_2(1):=b_2^*\L_2\ox\O_{B_2}(1)$.  Hence
$s^i(b_2^{-1}Z_1,b_2^*\L_1)B_2(\S)$ is a class on $b_2^{-1}Z_1\cap
B_2(S)$, and $b_{2*}s^j(\L_2(1))$ carries it to a class on $Z_1\cap S$,
so on $Z\cap S$.

View the joint blowup $\hat b_{1,2}\:\wh B_{1,2}\to\wh X$ as a repeated
blowup: the blowup $\hat b_2\:\wh B_2\to\wh X$ along $Z_2$, followed by
the blowup $q_1\:\wh B_{1,2}\to\wh B_2$ along $\hat b_2^{-1}Z_1$.  Then,
clearly,
	$$q_{1*}B_{1,2}(\S)=B_2(\S)^{X-Z_1}.$$
 Denote by $P^{2}_{X}$ and $P^{2}_{Z_2}$ the preimages in $\wh B_2$,
and by $P^{2}_{X}(\S)$ and $P^{2}_{Z_2}(\S)$ the corresponding cycles.
The projection formula yields $q_{1*}P^{1,2}_X(\S)=P^{2}_X(\S)$.
Therefore, \Cs3.1) and (4.3.1) yield
   $$q_{1*}P^{1,2}_Z(\S)=P^{2}_{Z_2}(\S)+B_{2}(\S)^{Z_1}.\tgs6.2$$

Rewrite \Cs4.1) using the identity,
  $$(\ell_1+\hat h_1)^i(\ell_2+\hat h_2)^j=\ell_1^i(\ell_2+\hat h_2)^j
  +(\ell_2+\hat h_2)^j\frac{(\ell_1+\hat h_1)^i-\ell_1^i}
	{\hat h_1}\hat h_1.$$
 Thus $s^{i,j}(Z_1,\L_1;Z_2,\L_2)(\S)$ is the sum of two terms.
Consider the first: it is equal, by the projection formula and \Cs6.2),
to
	$$\ell_1^{\,i}s^j(Z_2,\L_2)(\S)+
 b_{2*}\ell_1^{\,i}(\ell_2+ h_2)^jB_{2}(\S)^{Z_1}.\tgs6.3$$
 Consider the second term: it is equal, by  \Cs3.3), to
	$$b_*(\ell_2+ h_2)^j\frac{(\ell_1+ h_1)^i-\ell_1^i}
	{ h_1}D_1\cdot B_{1,2}(\S).\tgs6.4$$
 Now, $B_{1,2}$ is equal to the blowup of $B_2$ along $b_2^{-1}Z_1$.  So
(4.5.1) applies to $b_2^{-1}Z_1$ and $b_2^*\L_1$ and $B_2(\S)$; hence
\Cs6.4) is equal to
     $$b_{2*}(\ell_2+ h_2)^js^i(b_2^{-1}Z_1,b_2^*\L_1)B_{2}(\S)
    -b_{2*}\ell_1^{\,i}(\ell_2+ h_2)^jB_{2}(\S)^{Z_1}.\tgs6.5$$
 Finally,  \Cs6.3) and \Cs6.5) sum to the right side of \Cs6.1).
 \enddemo

\art 7 Relation to additivity

 The additivity formula (4.6.1) is a consequence of the mixed
operator formula  \Cs5.1) and the expansion formula \Cs6.1).  Indeed,
thanks to (4.4.3) (which is only an observation), they yield
 $$s^n(Z_1+Z_2,\L_1\otimes\L_2)=s^n(Z_2,\L_1\otimes\L_2)+
	   b_{2*}s^n(b_2^{-1}Z_1,\L_2(1)\otimes\L_1)B_{Z_2}. $$
 When $Z_2$ is a divisor, that formula is essentially (4.7.1), the
special case of (4.6.1) that is, as observed in (4.7), equivalent to
the general case.

Conversely, from the additivity formula (4.6.1), it is easy to
derive the mixed operator formula \Cs5.1) if the expansion formula
\Cs6.1) is assumed.

\sectionhead 6 Pseudo-divisors

\stp1
 On $X$, fix two invertible sheaves $\K$ and $\L$.  Fix a global
section $\kappa$ of $\K$, and denote its scheme of zeros by $K$.  The
triple $(\K,K,\kappa)$ is a `pseudo-divisor' (cf.~\cite{\Fulton,
\p.31}), and it can be intersected with any rational equivalence
class $\bold s$; however, we'll abuse notation by writing simply
$K\cdot\bold s$.

\lem2
 If $K$ contains no component of the cycle\/ $\S$ and no set
distinguished by $(Z,\S)$, then the three intersection cycles $
K\cdot\S$ and $K\cdot B_Z(\S)$ and $K\cdot P_Z(\S)$ are defined on $X$,
on $B_Z$ and on $\wh B_Z$ respectively, and they satisfy the following
two relations of commutativity:
 $$ K\cdot B_Z(\S)=B_Z( K\cdot\S)\and K\cdot P_Z(\S)=P_Z( K\cdot\S).
	\tgs2.1$$
 \pf
 It is easy to see that $ K\cdot\S$ and $K\cdot B_Z(\S)$ are defined
because $K$ contains no component of $\S$.  It is easy to see that
$K\cdot P_Z(\S)$ is defined because $K$ contains no set distinguished by
$(Z,\S)$.

Consider the first relation of \Cs2.1).  By construction, the right
hand cycle has no component contained in the exceptional divisor $D$.
Suppose that the left hand cycle has a component $\Gamma$ contained in
$D$.  By linearity, we may assume that $\S=[S]$ where $S$ is integral
with dimension $r$.  Then $\Gamma$ has dimension $r-1$, so is a
component of $D\cap S$.  Hence, by (3.7), the image $W$ of $\Gamma$ is
distinguished.  However, $W\subseteq K$.  No such $W$ can exist by
hypothesis.  Thus both sides of the first relation have no component
contained in $D$.  Hence, it suffices to prove that both sides agree on
$B_Z-D$.  Obviously, they do.

Consider the second relation.  To prove it, form the pullback $\wh K$
on $\wh X$.  Then
  $$\align
   K\cdot P_Z(\S)&=\wh K\cdot P_Z(\S)=\wh K\cdot P_Z \cdot B_Z(\wh\S\,)
			=P_Z\cdot\wh K\cdot B_Z(\wh\S\,)\\
		 &=P_Z\cdot B_Z(\wh K\cdot\wh\S\,)
	=P_Z\cdot B_Z\big(( K\cdot\S)\wh{\phantom K}\big)=P_Z( K\cdot\S).
 \endalign$$
  Indeed, the first equation holds because $ K=\wh K\cap X$, the second
and the last because $P_Z$ is a divisor, the third by commutativity of
intersection product, the fourth by the relation just proved, and the
fifth because $\wh X/X$ is flat.
 \enddemo

\propx{Commutativity with intersection product}3
 If $K$ contains no component of\/ $\S$ and no set
distinguished by $(Z,\S)$, then
	 $$ K\cdot s^i(Z,\L)(\S)=s^i(Z,\L)( K\cdot\S).$$
 \pf
 The left side is equal to $p_*(\ell+\hat h)^i K\cdot P_Z(\S)$ by (4.4.1)
and the projection formula.  So \Cs2) yields the assertion.
 \enddemo

\stp4
 Keep the setup of \Cs1).  In addition, fix $n\ge1$, and assume that
$\kappa$ arises from a global section of the subsheaf $\I^n\K$
where $\I$ is the ideal of $Z$.  Let $\kappa'$ be the section of
$\K(n)$ on $B_Z$ induced by the canonical surjection
$b^*\I\onto\O_B(1)$ where $b\:B_Z\to X$ is the blowup map, and let $
K'$ be its divisor of zeros.  Obviously, $ K\supseteq Z$ and
	$$b^* K = nD+ K' \tgs4.1$$
 where $D$ is the exceptional divisor.

\lem5
 Keep the setup of \Cs4).  Assume that $Z$ contains no component of
$\S$ and that $ K'$ contains no component of $B_Z(\S)$ or of $D\cdot
B_Z(\S)$.  Then the three intersection cycles $ K\cdot\S$ and $ K'\cdot
B_Z(\S)$ and $ K\cdot B_Z(\S)$ are defined, and the following relations
obtain:
	$$\gather
	 K'\cdot B_Z(\S)=B_Z( K\cdot\S)\tgs5.1\\
	 K\cdot B_Z(\S)=nD\cdot B_Z(\S)+ K'\cdot B_Z(\S)\tgs5.2\\
	b_*(nD\cdot B_Z(\S)) =( K\cdot \S)^Z\tgs5.3
	\endgather$$
 Moreover, every set distinguished by $(Z,\S)$ is also distinguished by
$(Z,\,K\cdot\S)$.
 \pf
 Obviously, $ K'\cdot B_Z(\S)$ is defined because $ K'$ contains no
component of $B_Z(\S)$.  Hence $ K\cdot B_Z(\S)$ is defined because of
\Cs4.1), and $ K\cdot\S$ is defined because $Z$ contains no component of
$\S$.

Consider the asserted relations.  First, $D$ contains no component of $
K'\cdot B_Z(\S)$ because $ K'$ contains no component of $D\cdot
B_Z(\S)$.  So it suffices to check \Cs5.1) off $D$, and there \Cs5.1) is
obvious because of \Cs4.1).  Second, \Cs4.1) yields \Cs5.2) trivially.
Third, \Cs5.2) and \Cs5.1) yield
  $$b_*(nD\cdot B_Z(\S)) =b_*( K\cdot B_Z(\S))-b_*(B_Z( K\cdot\S)).$$
 The third term is obviously equal to $( K\cdot\S)^{X-Z}$.  The second
term is, by the projection formula, equal to $ K\cdot b_*B_Z(\S)$, and
$b_*B_Z(\S)=\S$ because no component of $\S$ lies in $Z$.

To prove the last assertion, we may assume that $\S=[S]$ where $S$ is
integral with dimension $r$, and replace $X$, $Z$ and so forth by $S$,
$Z\cap S$ and so forth.  Let $W$ be a set distinguished by $(Z,\S)$.
Then $W$ is the image of a component $C$ of $P_Z$ with dimension $r$.
The hypothesis on $Z$ implies that $W$ is not a component of $X$.  So
$W$ is the image of a component of $D$ with dimension $r-1$ by (3.7).
Set $t:=\cod(W,X)$.  Then $t\ge1$ and $\dim W=r-t$ by (3.2.1).

If $t=1$, then $\dim W=r-1$, and so $W$ is a component of $Z$
distinguished by $(Z, K)$ by (3.3).  So assume $t\ge2$.  Then the
generic fiber $F$ of $C/W$ has dimension $t-1$ by \cite{\KT, Lemma
(3.2)(ii)}, and $t-1\ge1$.  Then $F\cap K'$ has pure dimension $t-2$
because it is the zero scheme of a section of the ample sheaf
$\O_F(n)$.  Consider a component $C'$ of the closure of $F\cap K'$.
Obviously, $C'$ maps onto $W$, and $\dim C'=r-2$ by \cite{\KT, Lemma
(3.2)(ii)}; hence, $C'$ is a component of $D\cap K'$.  Now, $ K'$ is
the proper transform of $ K$ because of \Cs4.1) and because $ K'$
contains no component of $D$ by hypothesis.  Hence (3.7) implies that
$W$ is distinguished by $(Z, K)$, and the proof is complete.
 \enddemo

\propx{Stability under intersection product}6
 Keep the setup of \Cs4).  Assume that $Z$  contains no
component of $\S$ and that $ K'$ contains no component of $B_Z(\S)$ or
of $D\cdot B_Z(\S)$.  Assume $\K=\L^{\ox n}$ and fix $j>0$.  Then
	$$ns^j(Z,\L)(\S)=s^{j-1}(Z,\L)( K\cdot\S).\tgs6.1$$
 Let $Z_1,\L_1$ be a second pair, and assume that $ K'$ contains no
set distinguished by $(b^{-1}Z_1,B_Z(\S))$.  Then
 $$ns^{i,j}(Z_1,\L_1;Z,\L)(\S)=s^{i,j-1}(Z_1,\L_1;Z,\L)( K\cdot\S).
 \tgs6.2
 $$
 \pf
 By hypothesis, $\S^Z=0$.  Hence, (4.5.2) yields
 $$ns^j(Z,\L)(\S)=
	\sum_{k+l=j-1}nb_*(\ell+h)^k\ell^lD\cdot B_Z(\S).$$
 By \Cs5.3) and the projection formula, $nb_*\ell^{j-1}D\cdot B_Z(\S)$ is
equal to $\ell^{j-1}( K\cdot\S)^Z$.  The remaining terms of the sum
have $k>0$, so contain $n(\ell+h)D\cdot B_Z(\S)$.  The later is, by
\Cs5.1), represented by the cycle $D\cdot B_Z(K\cdot\S)$.  Therefore,
again by (4.5.2), the sum is equal to $s^{j-1}(Z,\L)(K\cdot\S)$.  Thus
\Cs6.1) holds.

Consider \Cs6.2).  The left side is, by the expansion formula (5.6.1),
equal to
 $$n\ell_1^{\,i}s^j(Z,\L)(\S)+
	b_*n(\ell+h)^js^i(b^{-1}Z_1,b^*\L_1)B_Z(\S),\tgs6.3
 $$
 where $\ell_1:=c_1(\L_1)$.  By \Cs6.1), the first term is equal to
$\ell_1^is^{j-1}(Z,\L)( K\cdot\S)$.   By \Cs3) applied with
$B_Z(\S)$ and $\L(1)^{\ox n}$ for $\S$ and $\K$, the second term is
equal to
 $$b_*(\ell+h)^{j-1}s^i(b^{-1}Z_1,b^*\L_1)
  \bigl( K'\cdot B_Z(\S)\bigr).$$
 Again by \Cs5.1), that class is equal to the following one:
 $$b_*(\ell+h)^{j-1}s^i(b^{-1}Z_1,b^*\L_1)B_Z( K\cdot\S).$$
 Therefore, again by (5.6.1), the sum \Cs6.3) is equal to
the right side of \Cs6.2), as required.
 \enddemo

\sectionhead 7 Positivity

\lem1
 Assume $\S$ is an $r$-cycle, and let $W$ be a set
distinguished by $(Z,\S)$.  Set $t:=r-\dim W$.  Then the component
$s^i_W(Z,\L)(\S)$ in $A(W)$ of the twisted Segre class vanishes
for $i<t$, and, for $i\ge t$,
	$$\def\rng{k+l=i-1\atop l\ge t}
    s^i_W(Z,\L)(\S)=e\,\ell^{i-t}[W]+\sum_{\rng}\ell^k
	    b_*(\ell+h)^l(D\cdot B_Z(\S))\uprab W \tgs1.1$$
 for some integer $e$ independent of $i$, where $(D\cdot B_Z(\S))\uprab W $ is
the part of the intersection cycle whose components map onto $W$.

If $\S$ is prime, then $e>0$.  If $\S=[\Cal M]_r$ where $\Cal M$ is
coherent of dimension $ r$, then
	$$e=e(\I_{Z,w},\M_w) >0 \tgs1.2$$
 where $e(\I_{Z,w},\M_w)$ is the generalized Samuel multiplicity at
the generic point $w$ of $W$; if also $W$ is a component of
$Z\cap\Supp(\M)$, then $e$ is equal to the ordinary Samuel
multiplicity.
 \pf
 The proof of the blowup formula (4.5.2) yields
	$$s^i_W(Z,\L)(\S)=n_W\ell^i[W]
     +\sum_{k+l=i-1}b_*\ell^k(\ell+h)^l(D\cdot B_Z(\S))\uprab W \tgs1.3$$
 where $n_W$ is the multiplicity of $W$ as a component of $\S$.  If
$t=0$, then \Cs1.1) holds with $e=n_W$.  Assume $t>0$.  Then $n_W=0$.
Now, $(\ell+h)^l(D\cdot B_Z(\S))_W$ has dimension $r-1-l$.  Hence its
image under $b_*$ vanishes when $r-1-l>\dim W$, or equivalently, when
$l+1<t$.  Moreover, in $\A(W)$,
	$$b_*(\ell+h)^{t-1}(D\cdot B_Z(\S))\uprab W  = e[W]$$
 for some well-determined integer $e$.  Therefore, \Cs1.3) and the
projection formula yield \Cs1.1).

The value of $e$ may be found by making the flat base change to the
local scheme $X_w$.  Indeed, the formation of \Cs1.1) commutes with
this base change by standard theory; in particular,
	$$s^t_{\{w\}}(Z_w,\L_w)(\S_w)=e[w]$$
 where $Z_w$, $\L_w$ and $\S_w$ are the pullbacks to $X_w$ (note
that, in the case of an arbitrary Noetherian scheme $X$, when forming
the flat pullback $\S_w$, we must drop components that are of
dimension less than $t$).  Since $\L_w\simeq\O_{X,w}$, that equation
yields this one:
	$$\int{\hat h}^t[P_{Z_w}(\S_w)]^w=e$$
 where the superscript `$w$' means the part contained in the fiber
over $w$.  Suppose $\S=[\M]_r$.  Then $P_Z(\S)=[P_Z(\M)]_r$ by (4.2);
hence $[P_{Z_w}(\S_w)]_t$ is equal to $[P_{Z_w}(\M_w)]_t$.  Now,
$\dim\M_{w}=t$ and the ideal of $Z$ has maximal analytic spread on
$\M$ at $w$ by (3.2).  Hence the last two assertions hold by (3.6).
Finally, if $\S=[S]$ where $S$ is integral, then take $\M:=\O_S$, and
conclude $e>0$.
 \enddemo

\art 2 Multiplicity classes

Assume $\S$ is an $r$-cycle, and let $S$ be its support.  Then the
zero-dimensional rational equivalence class on $Z\cap S$ defined by
	$$e_k(\S)=e_k(Z,\L)(\S):=\ell^{k}s^{r-k}(Z,\L)(\S)$$
 will be called the $k$th  \dfn{BR-multiplicity class} of $\S$.
The similar class,
	$$m_k(\S)=m_k(Z,\L)(\S):=\ell^{k}s^{r-k}(Z,\O_X)(\S),$$
 will be called the $k$th {\it polar-multiplicity class}.  The two are
related by the following formulas, which follow directly from (4.4.3):
 $$e_k(\S)=\sum_{i=k}^{r}\binom{r-k}{r-i} m_{i}(\S)\and
	m_k(\S)=\sum_{i=k}^{r}(-1)^{i-k}\binom{r-k}{r-i}e_{i}(\S).$$

\prop3
 Assume $\S$ is an $r$-cycle.  Let $Y$ be a closed subscheme of $X$, and
$b_Y\:B_Y\to X$ and $b_Z\:B_Z\to X$ the blowup maps.  Then
  $$\align
  m_i(Y,\L)(\S)&=b_{Z*}m_i(b_Z^{-1}Y,\L_{B_Z}(1))B_Z(\S)\\
       &\qquad-b_{Y*}c_1(\O_{B_Y}(1))^{r-i}s^i(b_Y^{-1}Z,b_Y^*\L)B_Y(\S)
	\text{ for }i<r,\\
  m_r(Y,\L)(\S)&=e_0(Z,\L)(\S^Y)+b_{Z*}m_r(b_Z^{-1}Y,\L_{B_Z}(1))B_Z(\S).
	\endalign$$
 \pf
 By (5.4.2), $s^{i,j}(Z,\L;Y,{\Cal{O}}_X)(\S) =
s^{j,i}(Y,{\Cal{O}}_X;Z,\L)(\S)$.  So the expansion formula (5.6.1) yields
	$$\multline \ell^{i}s^j(Y)(\S)
	+b_{Y*}c_1(\O_{B_Y}(1))^js^i(b_Y^{-1}Z,b_Y^*\L)B_Y(\S) \\
	=s^j(\O_X)s^i(Z,\L)(\S)+b_{Z*}m_i(b_Z^{-1}Y,\L_{B_Z}(1))B_Z(\S).
	\endmultline$$
 Set $j:=r-i$.  For $i<r$, the assertion follows because $s^j(\O_X)$
vanishes.  For $i=r$, the assertion follows because of (4.4)(a) and
(4.4.4) and because
	$$e_0(Z,\L)(\S)-e_0(Z,\L)(\S^{X-Y})=e_0(Z,\L)(\S^Y).$$
 \enddemo

\art 4 Positive classes and positive sheaves

  Given an ambient scheme, call a rational equivalence class $\bold s$
\dfn{nonnegative\/} (resp., \dfn{positive\/}) and write $\bold s \scq
0$ (resp., $\bold s\succ0$) if some multiple $n\bold s$ with $n>0$ is
represented by a nonnegative cycle (resp., by a positive cycle).  A
positive cycle is nonzero by convention.  A positive class can vanish,
however, in some cases, but not if the ambient scheme is projective
over an Artin ring.

Call an invertible sheaf $\K$ \dfn{nonnegative\/} (resp.,
\dfn{positive\/}) and write $\K \scq 0$ (resp., $\K\succ0$) if
$c_1(\K)$ carries nonnegative (resp., positive) classes with
nonnegative (resp., positive) dimension into nonnegative (resp.,
positive) classes.  For example, $\K\scq0$ if $\K$ is generated by its
global sections; moreover, in the usual setup, $\L(1)|D\scq0$ if
$\I\L|Z$ is generated by its global sections where $\I$ is the ideal of
$Z$.

\prop5
 \text{\rm(1)}~~Assume $\L|Z\scq0$ and $\L(1)|D\scq0$.  If $\S\ge0$,
then
	$$s^i(Z,\L)(\S)\scq0\text{ for }i\ge0 \and
       e_0(Z,\L)(\S)\scq e_1(Z,\L)(\S)\scq\cdots\scq 0.\tgs5.1$$
 \part2 Assume $\L_1|Z\scq0$ where $Z:=Z_1+Z_2$, assume $\I_1\L_1$ is
generated along $Z_1$ by global sections, and assume $\L_2|Z_2\scq0$ and
$\L_2(1)|b_2^{-1}Z\scq0$.  If $\S\ge0$, then
	$$s^{i,j}(Z_1,\L_1;Z_2,\L_2)(\S)\scq0\text{ for }i,j\ge0.$$
  \pf
 To prove (1), set $\bold s^0:=\S^Z$ and $\bold s^{k+1}:=b_*
(\ell+h)^kD\cdot B_Z(\S)$ for $k\ge0$.  These classes are nonnegative
since $\L(1)|D\scq0$ and $\S\scq0$.  Now, the blowup formula
(4.5.2) and the projection formula yield
  $$s^i(Z,\L)(\S)=\sum_{0\le k\le i}\ell^{i-k}\bold s^k.$$
 Hence (1) follows since $\L|Z\scq0$.

Consider (2).  Part (1) yields $s^j(Z_2,\L_2)(\S)\scq0$ since
$\L_2|Z_2\scq0$ and $\L_2(1)|D_2\scq0$ as $D_2:=b_2^{-1}Z_2$.  Hence,
since also $\L_1|Z\scq0$, the first term in the expansion formula
(5.6.1) is nonnegative.  Its second term is also nonnegative, and to a
great extent the proof is similar.  Indeed, $\L_2(1)|b_2^{-1}Z_1\scq0$,
and $\L_1(1,0)\scq0$ on the blowup of $B_2$ along $b_2^{-1}Z_1$ (which
is equal to the joint blowup) since $\I_1\L_1$ is generated along $Z_1$
by global sections.  However, we don't assume $b_2^*\L_1|b_2^{-1}Z_1
\scq0$, but only $\L_1|Z_1\scq0$.  Nevertheless, the latter suffices
for carrying over, mutatis mutandis, the proof of (1).
 \enddemo

\prop6
 Assume $\L|Z\scq0$ and $\L(1)|D\scq0$.  Assume $\S$ is a positive
$r$-cycle, and let $W$ be a set distinguished by $(Z,\S)$.  Set
$d:=\dim W$ and $t:=r-d$.  Then
	$$s^t(Z,\L)(\S)\succ0.\tgs6.1$$
 Assume also $\L|Z\succ0$.  Then
 	$$\gather s^i(Z,\L)(\S)\succ0\text{ for }r\ge i\ge t,\tgs6.2\\
	e_k(Z,\L)(\S)\succ0\text{ for } k\le d.\tgs6.3\endgather$$
 Moreover, if either $\O_D(1)\scq0$ or $d$ is maximal, then
	$$m_d(Z,\L)(\S)\succ0.\tgs6.4$$
 \pf Let $i\ge t$.  Then \Cs1) yields $s^i_W(Z,\L)(\S)\scq
\ell^{i-t}[W]$.  Furthermore, clearly $s^i(Z,\L)(\S)\scq s^i_W(Z,\L)(\S)$.
Hence \Cs6.1) holds.  Assume $\L|Z\succ0$ also.  Then,
$\ell^{i-t}[W]\succ0$.  Hence \Cs6.2) holds.  Moreover,
$\ell^ks^{r-k}(Z,\L)(\S)\succ0$ for $r-k\ge t$, that is, for $k\le d$;
in other words, \Cs6.3) holds.  Finally, if also either $\O_D(1)\scq0$
or $d$ is maximal, then $s^t(Z,\O_X)(\S)\succ0$ either by \Cs6.1) or
because $s^t(Z,\O_X)(\S)\scq s^t_W(Z,\O_X)(\S)$; hence \Cs6.4) holds.
 \enddemo

\thm7
 Assume $\L_1|Z\succ0$ and $\L_1(1)|b_1^{-1}Z\scq0$.  Assume
$\L_2|Z\scq0$, and assume $\I_2\L_2$ is generated along $Z$ by global
sections.  Assume $\S$ is a positive $r$-cycle.  Let $W_2$ be a set
distinguished by $(Z_2,\S)$, and $W_1$ a set distinguished by
$(Z_1,W_2)$.  Set $t_k:=r-\dim W_k$ for each $k$.   Then
	$$s^{i,j}(Z_1,\L_1;Z_2,\L_2)(\S)\succ0
	 \text{ for }t_1\le i+j\le r\and j\le t_2.$$
 \pf
 Because of \Cs5)(2), we may assume, by linearity, that $\S=[S]$ where
$S$ is integral.  The proof proceeds by induction on $j$.  Suppose
$j=0$.  Now, (5.6.1) yields, thanks to (4.4.4) and (4.4)(a), this
equation:
  $$s^{i,0}(Z_1,\L_1;Z_2,\L_2)(\S)
	=s^i(\L_1)(\S^{Z_2})+s^i(Z_1,\L_1)(\S^{X-Z_2}).\tgs7.1$$
 If $t_2=0$, then $W_2=S$.  Since $W_2$ is contained in $Z_2$, the
first term in \Cs7.1) is, therefore, strictly positive.  So assume
$t_2>0$.  Then $S$ is not contained in $Z_2$, because, otherwise, $S$
would be the only distinguished subset of $Z_2$, but $W_2\ne S$.
Hence, on the right side of \Cs7.1), the first term vanishes, and the
second is equal to $s^i(Z_1,\L_1)(\S)$.  We'll prove in the next
paragraph that $W_1$ lies in a set $W$ distinguished by $(Z_1,\S)$.
Then \Cs6.2) will imply $s^i(Z_1,\L_1)(\S)\succ0$, completing the case
$j=0$.

To prove that $W$ exists, it suffices by (3.2) to prove that
	$$\dim W_1 + \cod(W_1,S)=r.\tgs7.2$$
 However, by (3.2),
   $$\dim W_2+\cod(W_2,S)=r\and\dim W_1+\cod(W_1,W_2)=\dim W_2.$$
 Hence
	$$r=\dim W_1+\cod(W_1,W_2)+\cod(W_2,S)
	\le\dim W_1+\cod(W_1,S)\le r.$$
 Thus \Cs7.2) holds, and the case $j=0$ is established.

Finally, assume $j>0$.  Then $t_2>0$.  Hence, as above, $S$ is not
contained in $Z_2$.  Therefore, by (6.6.2),
	$$ns^{i,j}(\S)=s^{i,j-1}( K\cdot\S)\tgs7.3$$
 where $n$ is a suitable integer and $ K$ is the scheme of zeros of a
suitable section $\kappa$ of $\L_2^{\ox n}$.  Such a suitable $n$ and
suitable section $\kappa$ exist because $\I_2\L_2$ is generated by its
global sections.  Indeed, the requirement on $n$ and $\kappa$ is
satisfied if $\kappa$ arises from a section of $\I_2^n\L_2^{\ox n}$ such
that the induced section of $\L_2^{\ox n}(n)$ does not vanish on certain
finitely many subsets of $D_2$; so a form of ``prime avoidance'' does the
trick.  Now, $W_2$ is also distinguished by $(Z_2, K\cdot\S)$ by the
last assertion of (6.5).  Therefore, \Cs7.3) and induction on $j$ yield
the assertion.
 \enddemo

\cor8
  Assume $\L_1|Z\succ0$ and $\L_1(1)|D_1\scq0$.  Assume $\I_2\L_2$ is
generated along $Z$ by global sections.  Assume $\S$ is a positive
$r$-cycle, and let $S$ be its support.  Assume that $Z_1\cap Z_2\cap S$
contains an irreducible closed set $W$ such that $\dim W+\cod(W,S)= r$.
Set $t_1:=r-\dim(W)$ and $t_2:=r-\dim(Z_2\cap S)$.  If
$\L_1(1)|b_1^{-1}Z\scq0$ and $\L_2|Z\scq0$, then
 	$$s^{i,j}(Z_1,\L_1;Z_2,\L_2)(\S)\succ0
	 \text{ for }t_1\le i+j\le r\and j\le t_2.$$
 \pf
 The assertion follows from \Cs7) since appropriate $W_k$ exist by
(3.2).
 \enddemo

\cor9
 Assume $\L|Z\succ0$ and assume $\L(1)|D\scq0$.  Let $Y$ be a second
closed subscheme of $X$, assume its ideal is generated along $Z$ by
global sections, and assume $\L(1)|b^{-1}Y\scq0$.  Let $\S$ be a
positive $r$-cycle.  Assume that there is a set $W$ distinguished by
$(Y,\S)$ and that there is some set distinguished by $(Z,W)$.  Set
$d:=\dim W$ and let $b_Z\:B_Z\to X$ be the blowup map.  Then
	$$b_{Z*}m_i(b_Z^{-1}Y,\L(1))B_Z(\S)\succ0\text{ for }r> i\ge d.$$
 \pf
 Theorem \Cs7) yields $s^{i,r-i}(Z,\L;Y,{\Cal{O}}_X)(\S)\succ 0$ for $r\ge
i\ge d$.  However, if $r> i$, then that Segre class is equal to the
class in question by the expansion formula (5.6.1) because
$s^{r-i}(\O_X)$ vanishes.
 \enddemo

\art 10 Relation to earlier work

The multiplicity classes $e_k(\S)$ and $m_k(\S)$ recover two of the
multiplicities of \cite{\KT}.  To establish the setup of the latter
paper in the present notation, take $X$ to be projective over a
Noetherian local ring, $\L$ to be $\O_X(1)$, and $Z$ to be defined by
forms of degree $1$.  In addition, let $Y$ be the preimage in $X$ of a
closed subscheme of the base supported at the closed point, and assume
that $Y$ contains the set $Z\cap S$ where $S$ is the support of $\S$.

Then the Buchsbaum-Rim multiplicity of \cite{\KT, (5.1)} is equal to the
degree of the zero-dimensional class $e_0(Z,\L)(\S)$ because of the
blowup formula (4.5.2) and because of \cite{\KT, (2.2.1)}.  In fact, for
every $n$, the $n$th associated multiplicity of \cite{\KT, (7.1)} is
equal to the degree of $e_n(Z,\L)(\S)$ by the same token.  Furthermore,
the additivity theorem (4.6) yields the additivity theorem \cite{\KT,
(6.7b)(i)}.  Similarly, the $n$th polar multiplicity of \cite{\KT,
(8.1)} is equal to the degree of $m_{r-n}(Z,\L)(\S)$.

The mixed multiplicities of \cite{\KT, (9.1)} are related to the mixed
Segre operators as follows.  Set
  $$x^{i,j}(\S):=b_{Y,*}h_Y^js^i(b_Y^{-1}Z,b_Y^*\L)B_Y(\S)\tgs10.1$$
  where $b_Y\:B_Y\to X$ is the blowup map and $h_Y:=c_1\O_{B_Y}(1)$.
Then (4.4)(a) yields
	$$x^{i,0}(\S)=s^i(Z,\L)(\S^{X-Y}).\tgs10.2$$
 Now, $B_Y(\S)^Z$ vanishes since $Z\cap S\subseteq Y$.  So the blowup
formula (4.5.1) yields
	$$x^{i,j}(\S)=b^{1,2}_*h_Y^j\frac{(\ell+h)^i-\ell^i}{h}D\cdot
		B_{1,2}(\S) \tgs10.3$$
 where $b^{1,2}\:B_{1,2}\to X$ is the joint blowup map of $Y$ and $Z$
and where $D$  is the preimage of $Z$.
Assume $i+j=r$.  Then $x^{i,j}(\S)$ has dimension zero, and as is
obvious from \Cs10.3), its degree is equal to the $j$th mixed
multiplicity of \cite{\KT, (9.1)}.  Moreover, the expansion formula
(5.6.1) yields
 $$s^{i,j}(Z,\L;Y,\O_X)(\S)=m_i(Y,\L)(\S)+x^{i,j}(\S).\tgs10.4
 $$
 Note that, by reason of dimension,
   $$m_i(Y,\L)(\S)=0\text{ for }i> d\text{ where }d:=\dim (Y\cap S).
	\tgs10.5$$

The mixed multiplicity $e^{i,j}(\S)$ of \cite{\KT, (9.10)(ii)}, where
$i+j=r$, is, thanks to (5.4.4), equal to the degree of the
zero-dimensional class $s^{i,j}(Z_1,\L_1;Z_2,\L_2)(\S)$.  So the
formula \Cs10.4) above, the symmetry relation (5.4.2), and the mixed
operator formula (5.5) yield the corresponding equations announced at
the end of \cite{\KT, (9.10)(ii)}.  Moreover, \Cs3) implies \cite{\KT,
(9.6)} because the latter's number $m^j(\bold T)$ is simply the degree
of $b_{Z*}m_i(b_Z^{-1}Y,(b_Z^*\L)(1))B_Z(\S)$.

The results of this section imply the main positivity results of
\cite{\KT}.  Indeed, since, by hypothesis, $Y$ contains the set $Z\cap
S$, we may replace $Z$ by the scheme $Z\cap Y$; then a positive
zero-cycle on $Z$ has positive degree, and so the present positivity
results yield positivity results about the multiplicities of
\cite{\KT}.  To be precise, assume $\S$ is positive.  If
$r=\dim\O_{S,y}$ for some (closed) $y\in Z\cap S$, then $y$ lies in a
set $W$ distinguished by $(Z,\S)$ by (3.2), and therefore
$e_0(Z,\L)(\S)\succ0$ by \Cs6).  Conversely, if $e_0(Z,\L)(\S)\succ0$,
then, obviously, there must exist a set distinguished by $(Z,\S)$, and
so, by (3.2), there exists a $y\in Z\cap S$ such that $r=\dim\O_{S,y}$.
Thus \cite{\KT, (5.2)(ii)} holds.  If $r=\dim\O_{S,y}$ for every $y\in
Z\cap S$, then every component of $Z\cap S$ is distinguished by
$(Z,\S)$ by (3.3), and therefore $e_d(Z,\L)(\S)\succ0$, where $d:=\dim
Z\cap S$ by \Cs6); thus the main assertion of \cite{\KT, (7.3)} holds.
If $r=\dim\O_{S,y}$ for every (closed) $y\in S$, then every component
of $Y\cap S$ is distinguished by $(Y,\S)$ by (3.3) since every closed
point of $X$ lies in $Y$, and therefore $m_d(Z,\L)(\S)\succ0$, where
$d:=\dim Y\cap S$, by \Cs6); thus \cite{\KT, (8.3)(vi)} holds.

Suppose again that $r=\dim\O_{S,y}$ for some (closed) $y\in Z\cap S$.
Now, $Z\cap S\subseteq Y$; hence, (3.2) implies that $y$ lies in a set
$W$ which is distinguished by $(Y,\S)$ and then that $y$ lies in a set
$W_1$ distinguished by $(Z,W)$.  (In fact, here any component of $Z\cap
W$ containing $y$ will serve as $W_1$ by (3.3).)  Hence \Cs10.4),
\Cs10.5) and \Cs7) imply $x^{i,r-i}(\S)\succ0$ for $r\ge i> d$; in
other words, \cite{\KT, (9.5)(3)$\Rightarrow$(1)} holds.  On the other
hand, \Cs9) yields $m^{r-i}(\bold T)>0$ for $r> i\ge d_W$ where
$d_W:=\dim W$; this conclusion improves \cite{\KT,
(9.7)(ii)(3)$\Rightarrow$(1)}, where the lower bound is $1+\dim(Y\cap
S)$, not $d_W$.  In particular, when $d<r$, we obtain $m^1(\T)>0$, the
inequality playing the key role in the proof of \cite{\KT, (10.1)}.

\sectionhead 8 Buchsbaum--Rim Polynomials

\stp1
 We'll use a mildly more general version of the setup of \cite{\KT}.
Thus $X=\Proj G$ where $G=\bigoplus G_n$ is a graded algebra over a
Noetherian local ring and $G$ is generated by finitely many elements of
$G_1$.  Then $Z$ is defined by (infinitely) many homogeneous ideals;
pick one, pick a system of generators, let $d$ be the maximum degree,
possibly 0, and let $H$ be the piece of degree $d$ of the ideal.  Then
the canonical map $H_X(-d)\to\O_X$ has as its image the (coherent)
ideal $\I$ of $Z$.

 Recall that the completed normal cone $P_Z$ arises from the graded
$\O_X$-algebra $\bigl(\bigoplus \I^p/\I^{p+1}\bigr)[u]$ where $u$ is an
indeterminate.  Twist this algebra by $\O_X(d)$; that is, tensor its
$p$th graded piece with $\O_X(pd)$.  The twist does not change the
scheme $P_Z$, but the new $\O_{P_Z}(1)$ is the old tensored with
$\O_X(d)$.  The $p$th power $\I^p$ arises from the homogeneous ideal of
$G$ generated by the $p$th power $H^p$.  Hence $P_Z$ arises from the
``twisted'' graded algebra whose $p$th graded piece is the $O_X$-module
arising from the graded $G$-module,
  $$\bigl(G/H^1G\oplus H^1G/H^2G\oplus\cdots H^pG/H^{p+1}G\bigr)[pd].
	\tgs1.1 $$
 Hence, $P_Z$ arises from the bigraded algebra, over the ground ring,
whose $(p,n)$-th bigraded piece is the $n$th graded piece of \Cs1.1),
namely,
	$$\bigoplus_{\nu=0}^p H^\nu G_{d(p-\nu)+n}
	/H^{\nu+1}G_{d(p-\nu-1)+n}.\tgs1.2$$
 Moreover, if the $\O_X$-module $\M$ arises from a graded $G$-module
$M$, then similarly the transform $P_Z(\M)$ arises from the bigraded
module whose $(p,n)$-th bigraded piece is
	$$\bigoplus_{\nu=0}^p H^\nu M_{d(p-\nu)+n}
	/H^{\nu+1}M_{d(p-\nu-1)+n}.\tgs1.3$$

 Let $\S$ be an $r$-cycle on $X$, and $S$ its support. Assume
$Z\cap S$ is contained in the closed fiber of $X$.  Then the twisted
Segre class $s^i(Z,\L)(\S)$ has support in the closed fiber, and so, if
$i+k=r$, then the following intersection number is defined:
	$$e^{i,k}(\S):=\int c_1\O_X(1)^ks^i(Z,\L)(\S).$$
 In the case where $\L=O_X(d)$, the number $e^{i,k}(\S)$ will be called
the $(i,k)$-th {\it Buchsbaum--Rim\/} multiplicity of $\S$.  Then the
zero-dimensional BR-multiplicity class $e_k(\S)$, defined in (7.2), is,
obviously, of degree $d^ke^{i,k}(\S)$.

\prop2
 Keep the setup of \Cs1).  Assume that the graded $G$-module $M$ is
finitely generated, that its associated sheaf $\M$ has dimension at most
$r$ and that the intersection $Z\cap\Supp\M$ is contained in the closed
fiber of $X$.  Then, as a function of $p$ and $n$, the length,
	$$\lambda(p,n):=\length(M_{pd+n}/H^pM_n),$$
 is eventually a polynomial of total degree at most $r$, and its term of
total degree $r$ is the form,
	$$\Lambda(p,n):=\sum_{i+k=r}
		e^{i,k}([\M]_r) p^i n^k/i!\,k!,$$
 where the coefficients  are the Buchsbaum--Rim multiplicities.
 \pf
 Use the new ``twisted'' $\O_{P_Z}(1)$.  Then the definitions yield
   $$e^{i,k}([\M]_r)=\int c_1\O_X(1)^k c_1\O_{P_Z}(1)^iP_Z([\M]_r).$$
 Now, $P_Z([\M]_r)=[P_Z(\M)]_r$ by (4.2), and $P_Z(\M)$ arises from the
bigraded module whose $(p,n)$-th bigraded piece is \Cs1.3).  Hence, by
the theory of Hilbert polynomials (see \cite{\KT, Lemma (4.3)} for
example), the length of \Cs1.3) is eventually a polynomial $\chi(p,n)$
of total degree at most $r$, and its term of total degree $r$ is the
form $\Lambda(p,n)$.  Obviously, the summands in \Cs1.3) are the
factors of a filtration of $M_{n+pd}/H^{p+1}M_{n-d}$.  So
$\lambda(p,n)$ is equal to $\chi (p-1,n+d)$.  However, these two
polynomials have the same term of total degree $r$.  Thus the assertion
holds.
 \enddemo

\stp3
 Keep the setup of \Cs1).  Assume that $Z$ is the sum $Z_1+Z_2$ in $X$
of two closed subschemes.  For each $Z_k$, pick a homogeneous ideal,
{}~pick a system of generators, let $d_k$ be the maximum degree, possibly
0, and let $H_k$ be the piece of degree $d_k$ of the ideal.  Recall
from (2.7) that $P^{1,2}_Z$ arises from a certain bigraded
$\O_X$-algebra.  As in \Cs1), it is convenient to twist this algebra;
this time, the $(p,q)$-th graded piece is tensored by
$\O_X(d_1p+d_2q)$.  The twist does not change the scheme $P^{1,2}_Z$,
but the first (resp., the second) tautological invertible sheaf is
replaced by its tensor product with $\O_X(d_1)$ (resp., $\O_X(d_2)$).

Recall from (2.8) that the transform $P^{1,2}_Z(\M)$ arises from a
certain bigraded $\O_X$-module, whose $(p,q)$-th graded piece is the
direct sum of the factors of certain finite filtration of
$\M/\I_1^{p+1}\I_2^{q+1}\M$.  As in \Cs1), it follows that $P^{1,2}_Z$
arises from a certain trigraded algebra over the ground ring, and that
$P^{1,2}_Z(\M)$ arises from a certain trigraded module, whose
$(p,q,n)$-th graded piece is the direct sum of the factors
of certain finite filtration of
	$$M_{d_1p+d_2q+n}/H^{p+1}_1H^{q+1}_2M_{n-d_1-d_2}$$
 where $M$ is any graded module giving arise to $\M$.

The mixed twisted Segre class has support in the closed fiber, and so,
if $i+j+k=r$, then the following intersection number is defined:
  $$e^{i,j,k}(\S):=\int c_1\O_X(1)^ks^{i,j}(Z_1,\L_1;Z_2,\L_2)(\S).$$
 It will be called the $(i,j,k)$-th {\it mixed Buchsbaum--Rim
multiplicity\/} of $\S$ when $\L_1=O_X(d_1)$ and $\L_2=O_X(d_2)$.  In
particular, $e^{i,j,0}(\S)$ is, as was noted in the middle of (7.10),
the mixed multiplicity $e^{i,j}(\S)$ of \cite{\KT, (9.10)(ii)}.

\thm4
 Keep the setup of \Cs3).  Assume that the graded $G$-module $M$ is
finitely generated, that its associated sheaf $\M$ has dimension at most
$r$ and that the intersection $Z\cap\Supp\M$ is contained in the closed
fiber of $X$.  Then, as a function of $p$, $q$  and $n$, the length,
	$$\lambda(p,q,n):=\length(M_{d_1p+d_2q+n}/H^{p}_1H^{q}_2M_n),$$
 is eventually a polynomial of total degree at most $r$, and its term of
total degree $r$ is the form,
	$$\Lambda(p,q,n):=\sum_{i+j+k=r}
		e^{i,j,k}([\M]_r) p^i q^j n^k/i!\,j!\,k!,$$
 where the coefficients  are the mixed Buchsbaum--Rim multiplicities.
 \pf
 The proof is entirely analogous to that of \Cs2).
 \enddemo

\art 5 Relation to earlier work

Continuing the discussion of (7.10), now to recover the results about
Buchsbaum--Rim polynomials in \cite{\KT}, keep the setup of \Cs1) and
take $\L$ to be $\O_X(1)$, and $Z$ to be defined by forms of degree $1$.
Then \Cs2) with $d:=1$ yields the conclusion of \cite{\KT, (5.8)(iii)};
in particular, with $n$ fixed suitably large, it yields the conclusion
of \cite{\KT, (5.7)}.  Let $Y$ be the preimage in $X$ of a closed
subscheme of the base defined by a primary ideal $J$.  Then \Cs2) with
$d:=0$ and $H:=J$ yields the main conclusion drawn in the middle of
\cite{\KT, (8.5)}.  Finally, \Cs4) with $d_1:=0$ and $H_1:=J$ and
$d_2=1$ and $H_2=H$ yields the main conclusions of \cite{\KT,
(9.10)(i)}.

\Refs

\atref{AM}
 \by R. Achilles and M. Manaresi
 \paper Multiplicity for ideals of maximal analytic spread and
intersection theory
 \jour J. Math. Kyoto Univ.\vol33--34 \yr1993 \pages1029--1046
 \endref

\atref{BR}
 \by D. A. Buchsbaum and D. S. Rim
 \paper A generalized Koszul complex. II. Depth and multiplicity
 \jour Trans. Amer. Math. Soc.\vol111 \yr1964 \pages197--224
 \endref

\atref{Fulton}
 \by W. Fulton
 \book Intersection Theory
 \bookinfo Ergebnisse der Mathematik und ihrer Grenzgebiete, 3.
Folge~$\cdot$ Band 2
 \publ Springer--Verlag\publaddr Berlin \yr1984
 \endref

\atref{FL}
 \by W. Fulton and D. Laksov
 \paper Residual intersections and the double point formula
 \inbook Real and complex singularities. Procedings, Oslo 1976
 \bookinfo P. Holm (ed.)
 \publ Sijthoff \& Noordhoff \yr1977 \pages171--77
 \endref

\atref{EGA}
 \by A. Grothendieck, with J.  Dieudonn\'e
 \book El\'ements de G\'eom\'etrie
 Alg\'ebrique II, III$_1$, IV$_1$, IV$_2$, IV$_3$
 \bookinfo Publ.  Math. I.H.E.S. {\bf 8}, {\bf 11}, {\bf 20}, {\bf 24},
  {\bf 26} \yr 1961, '61, '64, '65, '66
 \endref

\atref{KT}
\by S. Kleiman and A. Thorup
\paper A geometric theory of the Buchsbaum--Rim multiplicity
\jour J. Algebra \vol167\yr1994\pages168--231
\endref

\atref{NR}
 \by D. Northcott and D. Rees
 \paper Reductions of ideals in local rings
 \jour Math. Proc. Camb. Phil. Soc. \vol50 \yr1954 \pages145--58
 \endref

\atref{Teissier}
 \by B. Teissier
 \paper Cycles \'evanescents, sections planes et conditions de Whitney
 \paperinfo Singularit\'es \`a Car\-g\`ese
 \jour Ast\'erisque \vol7--8 \yr1973 \pages285--362
 \endref

\atref{Thorup}\by A. Thorup
\paper Rational equivalence on general noetherian schemes
\inbook Enumerative Geometry --  Sitges 1987
\bookinfo Proceedings, ed. S. Xambo-Descamps
\publ Springer-Verlag\publaddr New-York\yr1990\pages256--297
\endref

\endRefs
\enddocument